\documentclass{amsart} 
\usepackage{latexsym} 
\usepackage{amssymb} 
\usepackage{amsmath}
\usepackage{enumerate}
\usepackage{hyperref}  

\begin{document}


\newtheorem{theorem}{Theorem}[section]
\newtheorem{problem}{Problem} [section]
\newtheorem{definition}{Definition} [section]
\newtheorem{lemma}{Lemma}[section]
\newtheorem{proposition}{Proposition}[section]
\newtheorem{corollary}{Corollary}[section]
\newtheorem{example}{Example}[section]
\newtheorem{conjecture}{Conjecture} 
\newtheorem{algorithm}{Algorithm} 
\newtheorem{exercise}{Exercise}[section]
\newtheorem{remarkk}{Remark}[subsection]
 
\newcommand{\be}{\begin{equation}} 
\newcommand{\ee}{\end{equation}} 
\newcommand{\bea}{\begin{eqnarray}} 
\newcommand{\eea}{\end{eqnarray}} 

\newcommand{\eeq}{\end{equation}} 

\newcommand{\eeqn}{\end{eqnarray}} 
\newcommand{\beaa}{\begin{eqnarray*}} 
\newcommand{\eeaa}{\end{eqnarray*}} 

\newcommand{\lip}{\langle} 
\newcommand{\rip}{\rangle}

\newcommand{\uu}{\underline} 
\newcommand{\oo}{\overline} 
\newcommand{\La}{\Lambda} 
\newcommand{\la}{\lambda} 
\newcommand{\eps}{\varepsilon} 
\newcommand{\om}{\omega} 
\newcommand{\Om}{\Omega} 
\newcommand{\ga}{\gamma} 
\newcommand{\rrr}{{\Bigr )}} 
\newcommand{\qqq}{{\Bigl\|}} 
 
\newcommand{\dint}{\displaystyle\int} 
\newcommand{\dsum}{\displaystyle\sum} 
\newcommand{\dfr}{\displaystyle\frac} 
\newcommand{\bige}{\mbox{\Large\it e}} 
\newcommand{\integers}{{\Bbb Z}} 
\newcommand{\rationals}{{\Bbb Q}} 
\newcommand{\reals}{{\rm I\!R}} 
\newcommand{\realsd}{\reals^d} 
\newcommand{\realsn}{\reals^n} 
\newcommand{\NN}{{\rm I\!N}} 
\newcommand{\DD}{{\rm I\!D}} 
\newcommand{\degree}{{\scriptscriptstyle \circ }} 
\newcommand{\dfn}{\stackrel{\triangle}{=}} 
\def\complex{\mathop{\raise .45ex\hbox{${\bf\scriptstyle{|}}$} 
     \kern -0.40em {\rm \textstyle{C}}}\nolimits} 
\def\hilbert{\mathop{\raise .21ex\hbox{$\bigcirc$}}\kern -1.005em {\rm\textstyle{H}}} 
\newcommand{\RAISE}{{\:\raisebox{.6ex}{$\scriptstyle{>}$}\raisebox{-.3ex} 
           {$\scriptstyle{\!\!\!\!\!<}\:$}}} 
 
\newcommand{\hh}{{\:\raisebox{1.8ex}{$\scriptstyle{\degree}$}\raisebox{.0ex} 
           {$\textstyle{\!\!\!\! H}$}}} 

\newcommand{\OO}{\won} 
\newcommand{\calA}{{\mathcal A}} 
\newcommand{\calB}{{\cal B}} 
\newcommand{\calC}{{\cal C}} 
\newcommand{\calD}{{\cal D}} 
\newcommand{\calE}{{\cal E}} 
\newcommand{\calF}{{\mathcal F}} 
\newcommand{\calG}{{\cal G}} 
\newcommand{\calH}{{\cal H}} 
\newcommand{\calK}{{\cal K}} 
\newcommand{\calL}{{\mathcal L}} 
\newcommand{\calM}{{\mathcal M}} 
\newcommand{\calO}{{\cal O}} 
\newcommand{\calP}{{\cal P}} 
\newcommand{\calU}{{\mathcal U}} 
\newcommand{\calX}{{\cal X}} 
\newcommand{\calXX}{{\cal X\mbox{\raisebox{.3ex}{$\!\!\!\!\!-$}}}} 
\newcommand{\calXXX}{{\cal X\!\!\!\!\!-}} 
\newcommand{\gi}{{\raisebox{.0ex}{$\scriptscriptstyle{\cal X}$} 
\raisebox{.1ex} {$\scriptstyle{\!\!\!\!-}\:$}}} 
\newcommand{\intsim}{\int_0^1\!\!\!\!\!\!\!\!\!\sim} 
\newcommand{\intsimt}{\int_0^t\!\!\!\!\!\!\!\!\!\sim} 
\newcommand{\pp}{{\partial}} 
\newcommand{\al}{{\alpha}} 
\newcommand{\sB}{{\cal B}} 
\newcommand{\sL}{{\cal L}} 
\newcommand{\sF}{{\cal F}} 
\newcommand{\sE}{{\cal E}} 
\newcommand{\sX}{{\cal X}} 
\newcommand{\R}{{\rm I\!R}} 
\renewcommand{\L}{{\rm I\!L}} 
\newcommand{\vp}{\varphi} 
\newcommand{\N}{{\rm I\!N}} 
\def\ooo{\lip} 
\def\ccc{\rip} 
\newcommand{\ot}{\hat\otimes} 
\newcommand{\rP}{{\Bbb P}} 
\newcommand{\bfcdot}{{\mbox{\boldmath$\cdot$}}} 
 
\renewcommand{\varrho}{{\ell}} 
\newcommand{\dett}{{\textstyle{\det_2}}} 
\newcommand{\sign}{{\mbox{\rm sign}}} 
\newcommand{\TE}{{\rm TE}} 
\newcommand{\TA}{{\rm TA}} 
\newcommand{\E}{{\rm E\, }} 
\newcommand{\won}{{\mbox{\bf 1}}} 
\newcommand{\Lebn}{{\rm Leb}_n} 
\newcommand{\Prob}{{\rm Prob\, }} 
\newcommand{\sinc}{{\rm sinc\, }} 
\newcommand{\ctg}{{\rm ctg\, }} 
\newcommand{\loc}{{\rm loc}} 
\newcommand{\trace}{{\, \, \rm trace\, \, }} 
\newcommand{\Dom}{{\rm Dom}} 
\newcommand{\ifff}{\mbox{\ if and only if\ }} 
\newcommand{\nproof}{\noindent {\bf Proof:\ }} 
\newcommand{\nproofYWN}{\noindent {\bf Proof of Theorem~\ref{YWN}:\ }} 
\newcommand{\remark}{\noindent {\bf Remark:\ }} 
\newcommand{\remarks}{\noindent {\bf Remarks:\ }} 
\newcommand{\note}{\noindent {\bf Note:\ }} 
 \newcommand{\examples}{\noindent {\bf Examples:\ }} 
 
\newcommand{\boldx}{{\bf x}} 
\newcommand{\boldX}{{\bf X}} 
\newcommand{\boldy}{{\bf y}} 
\newcommand{\boldR}{{\bf R}} 
\newcommand{\uux}{\uu{x}} 
\newcommand{\uuY}{\uu{Y}} 
 
\newcommand{\limn}{\lim_{n \rightarrow \infty}} 
\newcommand{\limN}{\lim_{N \rightarrow \infty}} 
\newcommand{\limr}{\lim_{r \rightarrow \infty}} 
\newcommand{\limd}{\lim_{\delta \rightarrow \infty}} 
\newcommand{\limM}{\lim_{M \rightarrow \infty}} 
\newcommand{\limsupn}{\limsup_{n \rightarrow \infty}} 
 
\newcommand{\ra}{ \rightarrow } 

 \newcommand{\mlim}{\lim_{m \rightarrow \infty}}  
 \newcommand{\limm}{\lim_{m \rightarrow \infty}}  
 \newcommand{\nlim}{\lim_{n \rightarrow \infty}} 
 
 
 
 
 
 
 
\newcommand{\one}{\frac{1}{n}\:} 
\newcommand{\half}{\frac{1}{2}\:} 
 
\def\le{\leq} 
\def\ge{\geq} 
\def\lt{<} 
\def\gt{>} 
 
\def\squarebox#1{\hbox to #1{\hfill\vbox to #1{\vfill}}} 
\newcommand{\nqed}{\hspace*{\fill} 
           \vbox{\hrule\hbox{\vrule\squarebox{.667em}\vrule}\hrule}\bigskip} 

\title[Weak calculus of variations]{Weak  calculus of variations for functionals of laws of semi-martingales} 
\begin{abstract}
We develop a non-anticipating calculus of variations for functionals on a space of laws of continuous semi$-$martingales, which extends the classical one. We extend Hamilton's least action principle and Noether's theorem to this generalized stochastic framework. As an application we obtain, under mild conditions, a stochastic Euler$-$Lagrange condition and invariants for the critical points of recent problems in stochastic control, namely for the semi-martingale optimal transportation problems.   \end{abstract}
 
 \author{R\'{e}mi Lassalle}
\author{Ana Bela Cruzeiro}  

\maketitle 
\noindent

\vspace{0.5cm}

\noindent 

 \section*{Introduction}

In this paper we formulate a weak calculus of variations which extends the classical one.  Roughly speaking this enables to perform  a calculus on functions defined on laws of semi-martingales. We apply this calculus to obtain a stochastic extension of Hamilton's least action principle. Recall that the classical version of this principle provides a characterization of the paths satisfying the Euler$-$Lagrange condition as critical points of a functional which is called an action. Here we will characterize laws of semi-martingales which satisfy a constraint that extends the classical one.  Namely these laws will be proved to be critical points of a stochastic action. Once this extension is achieved we use it to relate some invariance  properties of the critical processes to the symmetries of the corresponding Lagrangian; in other words, we derive a stochastic extension of Noether's theorem. Finally we consider applications to stochastic control, in particular to some semi-martingale optimal transportation problems. These problems were recently introduced in  \cite{TOUZI} with application to financial mathematics. 

As a warm up, let us give, in an informal way, some details on our motivation and on the difficulties we overcome with our approach. The first motivation lies in classical mechanics. In classical mechanics one usually considers paths sufficiently regular to model the kinematics of a system. In particular if one  describes the trajectory of  a classical particle by  a path $q : [0,1] \to \mathbb{R}$ one will usually ask $q$ to be sufficiently regular in order to provide a realistic description of the observation. Namely it will be often assumed to be $C^2$ for both its speed $\dot{q}_t$ ($:= \frac{dq_t}{dt}$) and its acceleration $\ddot{q}_t$ to be defined. Thus, for the sake of simplicity let us first consider the space $\Omega^2_{[0,1]}$ of the $C^2$ paths $q : [0,1] \to \mathbb{R}$ as being the set of the paths providing an admissible description of the kinematics.  The possibility to make predictions i.e., to be able to estimate the configuration of the system $(q_t, \dot{q}_t)$  at time $t$ from the initial conditions, relies on the existence of a dynamics which is of physical origin. This latter is expressed in the model by a further constraint on the $q$ paths which involves a function  $$\mathcal{L} :(x,v) \in \mathbb{R}\times \mathbb{R} \to \mathcal{L}(x,v)\in \mathbb{R}$$ where  $x$ (resp. $v$) may stand for the position (resp. the speed). This function $\mathcal{L}$, which is called a Lagrangian, contains all the physics of the model, and the related constraint  which is called the Euler$-$Lagrange condition (see \cite{AM}, \cite{ARNOLD6}, \cite{LANDAU})  reads \begin{equation} \label{elcint1} \frac{d}{dt}\partial_v \mathcal{L}(q_t,\dot{q}_t)= \partial_x \mathcal{L}(q_t,\dot{q}_t) \end{equation}   Integrating in time, it becomes \begin{equation}  \label{stoemp0}\partial_v \mathcal{L}(q_t, \dot{q}_t) - \int_0^t \partial_x \mathcal{L}(q_s, \dot{q}_s) ds = c\end{equation} where $c$ is some constant. Under mild conditions on $\mathcal{L}$ the paths $q\in \Omega_{[0,1]}^2$ satisfying the Euler$-$Lagrange condition can be characterized as critical points of a functional $\mathcal{S}_{path}$ which is called the action of the system (see \cite{AM}). It is defined by $$\mathcal{S}_{path}(q) = \int_0^1 \mathcal{L}(q_t, \dot{q}_t) dt$$ and $q$ is said to be critical if for all $h\in  \Omega_{[0,1]}^2$ satisfying $h_0=h_1=0$  $$\frac{d}{d\epsilon }\mathcal{S}_{path}(q^\epsilon)|_{\epsilon=0}= 0$$ where for $\epsilon \in \mathbb{R}$, $q^\epsilon := q + \epsilon h$  is a perturbation of the path. The theorem which states the equivalence for a path $q\in \Omega^2_{[0,1]}$ to satisfy the Euler$-$Lagrange condition~(\ref{elcint1}) or to be a critical point of the action is called Hamilton's least action principle (see \cite{AM}, \cite{ARNOLD6},  \cite{LANDAU}). One of the goals of this paper is to extend this result to some stochastic framework.  

Let us denote by $\mathbb{S}$ the set of  laws of continuous semi-martingales such that for $\nu\in \mathbb{S}$, the canonical process satisfies $\nu-a.s.$ for all $t\in[0,1]$,
$$W_t = W_0 + \int_0^t v_s^\nu ds + M_t^\nu$$ where $(M_t^\nu)$ is some (local) martingale on the probability space $(C([0,1], \mathbb{R}),\mathcal{B}( C([0,1], \mathbb{R}))^\nu, \nu)$ for the filtration $(\mathcal{F}_t^\nu)$ (which denotes the $\nu-$usual augmentation of the filtration generated by the evaluation process),  where  $(<M^\nu >_t)$ is assumed to be absolutely continuous with a derivative ($\alpha_t^\nu$). Setting $\nu = \delta^{Dirac}_q$ for $q\in \Omega^2_{[0,1]}$ we have $\nu-a.s.$ for all $t\in[0,1]$  \begin{equation} \label{kinem} W_t = W_0 + \int_0^t v_s^\nu ds  \end{equation} where $\lambda\otimes \nu-a.s.$ $$v_t^\nu = \dot{q}_t$$  i.e. $\nu \in \mathbb{S}$ and $M^\nu = 0$.  Thus, let us regard $\mathbb{S}$ as an extension of the set of the paths describing admissible kinematics in an extended stochastic context. In this paper we will consider a constraint, which we call the stochastic Euler$-$Lagrange condition, that extends  on $\mathbb{S}$ the classical Euler$-$Lagrange condition ; in particular it is a natural way to introduce some dynamics in a stochastic framework (see also \cite{RLZ}).  Namely, given some suitable function $\mathcal{L} : (t, x,y,a) \in[0,1]\times \mathbb{R}\times \mathbb{R}\times \mathbb{R} \to \mathcal{L}_t(x,y,a)\in \mathbb{R}$ a law $\nu \in \mathbb{S}$ will  be said to satisfy the stochastic Euler$-$Lagrange condition if   \begin{equation}  \label{stoemp}\partial_v \mathcal{L}_t(W_t, v_t^\nu,\alpha_t^\nu) - \int_0^t \partial_x \mathcal{L}_s(W_s, v_s^\nu, \alpha_s^\nu) ds = N_t^\nu\end{equation}  for some $(\mathcal{F}_t^\nu)-$martingale $(N_t^\nu)$ on  $(C([0,1], \mathbb{R}),\mathcal{B}( C([0,1], \mathbb{R}))^\nu, \nu)$. Indeed by taking $\nu = \delta^{Dirac}_q$ for $q\in \Omega_{[0,1]}^2$ and a Lagrangian $\mathcal{L}$ not depending on $a$ and $t$~(\ref{stoemp}) is equivalent to~(\ref{stoemp0}). By extending Hamilton's least action principle to $\mathbb{S}$ we will relate the dynamical condition~(\ref{stoemp}) to recent problems of stochastic control which is our second motivation. 

Consider the variational problems of the form \begin{equation} \label{touzpb} \inf\left(\left\{ \mathcal{S}(\nu) : \nu \in \mathbb{S}, Law(W_0)= \nu_0, Law(W_1)= \nu_1\right\}\right) \end{equation}
where \begin{equation} \label{actiontouz} \mathcal{S}(\nu):= E_\nu\left[\int_0^1 \mathcal{L}(W_s,v_s^\nu, \alpha_s^\nu) ds \right]. \end{equation} 
Such problems (extending those considered in \cite{Mika1}, \cite{Mika3}) have been recently investigated in \cite{TOUZI}; one minimizes among laws of semi-martingales with fixed initial (resp. final) marginal law $\nu_0$ (resp. $\nu_1$). As a matter of fact they extend the so-called Sch\"{o}dinger problem (see \cite{F3} and \cite{LEO}), which can be written as an entropy minimization problem. In this latter case, where the optimal processes may be computed explicitly, it was noticed by J.C. Zambrini (see \cite{RLZ} for instance) that the optimum solves a stochastic Euler$-$Lagrange condition (\ref{stoemp}). On the other hand in the general case of~(\ref{touzpb}), or by considering even more general problems where one fixes the joint law  (see \cite{leo} for the case of Bernstein's processes) of $(W_0,W_1)$ to be equal to a given Borel probability $\gamma$ on $\mathbb{R}^d\times \mathbb{R}^d$, \begin{equation} \label{touzleopb} \inf\left(\left\{ E_\nu\left[\int_0^1 \mathcal{L}(W_s,v_s^\nu, \alpha_s^\nu) ds \right] : \nu \in \mathbb{S}, Law(W_0,W_1)= \gamma \right\}\right) \end{equation} It is not convenient to use explicit formulae: in this paper  we rather state a stochastic least action principle which extends the classical one, proving that the optimum of these problems of stochastic control are  critical points of a stochastic action. In the classical Hamilton's principle the paths satisfying  Euler$-$Lagrange conditions are critical points and not necessarily minimum. Similarly, within our stochastic extension we also allow processes satisfying~(\ref{stoemp}) which are not minimum for problems of the form~(\ref{touzleopb}).  Actually, as it will be pointed out on examples on the classical Wiener space for a quadratic cost,  the situation is more complicated in the stochastic case of~(\ref{stoemp}). 

 We then prove a Noether theorem, which we apply to the extremum of~(\ref{touzpb}) and~(\ref{touzleopb}).

 We found inspiration for applications to stochastic control essentially in \cite{ZT},\cite{Zamb1},  where they focus on Bernstein's processes. Our results may be compared to those.  We also show that in some cases~(\ref{stoemp}) is related to  systems of coupled stochastic differential equation and to PDEs (such as Navier$-$Stokes equations).

Finally, let us  add some comments concerning technical issues.  When one expresses the proof of the least action principle using probabilities by  \begin{equation} \label{wmple} \Omega_{[0,1]}^2 \hookrightarrow_{\delta^{Dirac}} \mathbb{S}\end{equation}  we set $\nu_\epsilon = \delta^{Dirac}_{q^\epsilon}$  and differentiate $$\frac{d}{d\epsilon}\mathcal{S}(\nu_\epsilon)|_{\epsilon=0}.$$ where $$\nu_\epsilon = (I_W + \epsilon h)_\star \nu$$ i.e. the variation becomes the image of $\nu$ by the measurable mapping $\tau_{\epsilon h} : \omega \in W \to \omega + \epsilon h$ for $h\in \Omega_{[0,1]}^2$. In a stochastic framework one will have to consider more general perturbations of the form $$\tau_k : \omega \in C([0,1], \mathbb{R}) \to \omega + k(\omega)\in C([0,1], \mathbb{R})$$ where $k:=\int_0^. \dot{k}_s ds$ is now random and adapted to the canonical filtration. Setting  \begin{equation} \label{nuemf} \nu_\epsilon := (I_W+ \epsilon k)_\star\nu \end{equation} we realize that some essential properties will not necessarily hold. We do not have that $\tau_k$ is invertible (even almost surely) in general, and most of all in general we do not have  a.e.  $$v_t^{\nu_\epsilon} (\omega+ \epsilon k(\omega))  = v_t^{\nu}(\omega) + \epsilon \dot{k}_t(\omega).$$  As a consequence we cannot differentiate relevant functionals in all  (adapted) random directions.  This is essentially due to the fact that such perturbations may not preserve the filtration. To overcome these difficulties we build, for any $\nu \in \mathbb{S}$, some associated vector space of variation processes, which is roughly speaking the set  directions towards which the variations  of relevant functionals on $\mathbb{S}$ can be handled as in the classical case. Then we prove that the space is wide enough to build a derivative on $\mathbb{S}$ and to obtain a necessary and sufficient condition for~(\ref{stoemp}) on $\mathbb{S}$ by means  of a least action principle.

The structure of the paper is the following. In Section~\ref{section1} we fix the notations of the whole paper and we recall the variational characterization of martingales, as well as some results about transformation of measures preserving the filtration. In the following section we define the variation processes and state their main properties, 
namely that they form a dense vector subspace of the space of the adapted shifts of finite energy. In Section~\ref{section3} we compute the changing formula of the characteristics 
of a $\nu\in \mathbb{S}$ given several particular transformations of measure  (which will be used to  compute explicitly the differential of actions on $\mathbb{S}$). We also lift transformations of  space depending on the time to transformations on $\mathbb{S}$.  In Section~\ref{section4} we define the  differential of functionals defined on $\mathbb{S}$ in such a way that extends the usual calculus of variations by~(\ref{wmple}). We note  that the definition extends directly to  Borel probabilities 
on the space of  continuous functions. In Section~\ref{section5} we state precisely the definition of the laws satisfying the stochastic Euler-Lagrange condition,  our hypothesis on  the Lagrangian (we call it regular) and we prove the stochastic least action principle Theorem~\ref{thmlap}, which is our main result. Then, in Section~\ref{section6}, we generalize
Noether's theorem (such as it is formulated by \cite{ARNOLD6}) to this general framework (Theorem~\ref{noether}). Section~\ref{section7} is devoted to applications in  stochastic 
control and in particular to the problems considered in \cite{TOUZI}. Namely we obtain some information on the optimum of variational problems by using the stochastic least action principle 
and Noether's theorem. Finally in the last section we illustrate the content of~(\ref{stoemp}) in the case  of the classical action defined on $\mathbb{S}$ and we investigate the corresponding critical processes. In this case we relate the results to systems of stochastic differential equations and provide some explicit examples and counterexamples.

\section{Preliminaries and notations}
\label{section1}

\subsection{The path spaces and their stochastic counterparts.}
\label{morph}
In the whole paper $(\Omega,\mathcal{A},\mathcal{P})$  will always denote a complete probability space  and $(\mathcal{A}_t)$ a filtration on $\Omega$ satisfying the usual conditions  (i.e. right continuous and complete) such that for all $t\in[0,1]$, $\mathcal{A}_t\subset \mathcal{A}$.  Under these hypothesis, following \cite{JACOD}, we call $(\Omega, \mathcal{A}, (\mathcal{A}_t)_{t\in[0,1]}, \mathcal{P})$ a complete stochastic basis. We emphasize that all these  assumptions are crucial for our results to hold. The most convenient way to handle transformations of  laws of stochastic processes whose trajectories are sufficiently regular is to consider them as random trajectories.   Thus consider the space  $W= C([0,1],\mathbb{R}^d)$ of continuous functions on $[0,1]$ with values in $\mathbb{R}^d$.  Processes will be often regarded as random elements taking their values  in $W$, and we will sometimes call the elements of $W$ paths or trajectories. 

We recall that $W$ is  a separable Banach space  with respect to the norm $|.|_W$ of the uniform  convergence ($|\omega|_W:= \sup_{t\in[0,1]}|\omega(t)|_{\mathbb{R}^d}$).  We can consider the related Borel sigma field $\mathcal{B}(W)$, which turns $W$ into a measurable space. Within this  perspective,  we consider a continuous stochastic process $(X_t)_{t\in[0,1]}$ as a $\mathcal{A}/ \mathcal{B}(W)-$measurable mapping  $X : \Omega \to W$.  

 We denote by $\mathcal{P}_W$  the set of Borel probabilities on $W$, which are the laws of the continuous processes seen as random trajectories. We denote  $f_\star \mathcal{P}$ the image of a measure $\mathcal{P}$ by a measurable mapping $f: \Omega \to \widetilde{\Omega}$ where $(\widetilde{\Omega}, \widetilde{\mathcal{A}})$ is some other measurable space.

 In the sequel we shall  work under the usual conditions that insure existence of sufficiently regular modifications of martingales. Therefore we will always work on complete probability spaces with filtrations satisfying the usual conditions (i.e. complete and right continuous).  Taking this into account we introduce the following notations. If $\eta \in \mathcal{P}_W$ and $\mathcal{G}$ is a sigma-field  such that $\mathcal{G}\subset \mathcal{B}(W)$  $\mathcal{G}^\eta$ will denote the $\eta-$completion of $\mathcal{G}$ i.e. the smallest sigma field which contains all the elements of $\mathcal{G}$ and all the $\eta-$negligible sets. The unique extension of $\eta \in \mathcal{P}_W$ to $\mathcal{B}(W)^\eta$ will be still denoted by $\eta$.   We denote by $(W_t)$ the evaluation process defined by  $$W_t : \omega \in W \to W_t(\omega)\in \mathbb{R}^d$$ for $t\in [0,1]$.  For $\eta \in \mathcal{P}_W$, $(W_t)$ defines a process on the  probability space $(W,\mathcal{B}(W)^\eta,\eta)$ :  it is how we will consider it in the sequel. The corresponding measurable mapping is the identity  $$I_W : \omega \in W \to \omega \in W$$ which is Borel measurable (and thus $\mathcal{B}(W)^\eta/ \mathcal{B}(W)-$measurable).  
 
 By considering a  path $\omega \in W$, and denoting by $\delta^{Dirac}_\omega\in \mathcal{P}_W$ the Dirac measure concentrated on $\omega$ (i.e. $\delta^{Dirac}_\omega(A)= I_A(\omega)$, $A\in \mathcal{B}(W)$) we obtain an embedding $$W\hookrightarrow_{\delta^{Dirac}} \mathcal{P}_W .$$ In this sense any path can be seen as a stochastic process, and the weak calculus of variations we will introduce below is such that, through this embedding, it extends the classical one. More generally transformations of measures can be  formalized by transference plans (Borel probabilities of $W\times W$). In this work we shall not need this generality: transformations of measure will be merely achieved by images of probabilities induced by measurable mappings. More precisely we will handle equivalence classes of mappings. For  $(\Omega, \mathcal{A}, \mathcal{P})$ a complete probability space, $M_{\mathcal{P}}((\Omega, \mathcal{A}),(W,\mathcal{B}(W))$  will denote the set  obtained by identifying  $\mathcal{A}/ \mathcal{B}(W)-$measurable mappings $f : \Omega \to W$ which are $\mathcal{P}-a.s.$ equal. Following \cite{MAL2} we will sometimes call the elements of this space morphisms of probability spaces. If $U\in M_{\mathcal{P}}((\Omega, \mathcal{A}),(W,\mathcal{B}(W))$ and $f : \Omega \to W$ is a $\mathcal{A}/ \mathcal{B}(W)-$measurable mapping we will note $\mathcal{P}-a.s.$ $U=f$ to denote that the $\mathcal{P}-$equivalence class associated to $f$ is $U$ (i.e. the $\mathcal{P}-$equivalence class $U$ can be seen as the set of the $\mathcal{A}/ \mathcal{B}(W)-$measurable mappings $g:\Omega\to W$ such that $\mathcal{P}-a.s. $ $f=g$). Similarly if $V\in M_{\mathcal{P}}((\Omega, \mathcal{A}),(W,\mathcal{B}(W))$ we will note $\mathcal{P}-a.s.$ $U=V$ to denote that $U$ and $V$ are the same $\mathcal{P}-$equivalence class.



%



 We introduce  the Hilbert space of the absolutely continuous functions on $[0,1]$ vanishing at $t=0$ with a square integrable derivative $$H:=\left\{ h : [0,1] \to \mathbb{R}^d , h:=\int_0^. \dot{h}_s ds , \int_0^1 |\dot{h}_s|_{\mathbb{R}^d}^2 ds <\infty\right\}$$ (the so-called Cameron-Martin space) and we note  $<.,.>_H$ (resp. $|.|_H$) the corresponding Hilbert product (resp. norm).  Then we denote by $W_{abs}$  the subset of $W$ whose elements are absolutely continuous functions (i.e. the set of $\omega\in W$ of the form $\omega := \int_0^.\dot{\omega}_s ds$) and by $H_{0,0}$ the subset of $H$ given by \begin{equation} \label{h00} H_{0,0}:= \left\{ h \in H : h_1=0 \right\} .  \end{equation} Note that by definition of $H$ for $h\in H_{0,0}$ we have $h_0=h_1=0$. In the classical setting this set is the set of variations. Our main task will be to build its counterpart in the stochastic framework, and we will need to consider spaces of (equivalence classes) of mappings taking almost surely their values in such spaces.

When $E$ is a Borel measurable subset of $W$, let us denote by $L^0(\mathcal{P}, E)$ the space of the $\mathcal{P}-$equivalence classes of 
mappings $u$ (i.e. $u\in M_{\mathcal{P}}((\Omega, \mathcal{A}),(W, \mathcal{B}(W)))$ such that $\mathcal{P}-a.s.$ $u \in E$. To control the 
integrability within this stochastic context, we also need the space $L^2(\mathcal{P}, H)$ (resp. $L^2(\mathcal{P}, H_{0,0})$) of the functions $u\in L^0
(\mathcal{P}, H)$ (resp. in $L^0(\mathcal{P}, H_{0,0})$) such that $E_{\mathcal{P}}[|u|_H^2] <\infty$ i.e. $$E_{\mathcal{P}}\left[ \int_0^1 |\dot{u}_s|_{\mathbb{R}^d}^2 ds \right] < \infty$$ where $\mathcal{P}-a.s.$ $u=\int_0^. \dot{u}_s ds$. Similarly  $L^\infty(\mathcal{P}, W)$ (resp. $L^\infty(\mathcal{P},H)$) will denote the set of  $u\in L^0(\mathcal{P},W)$ for which there exists a $K_u>0$ such that $\mathcal{P}-a.s.$ $|u|_W<K_u$ (resp. $|u|_H<K_u$). One of the main differences with respect to the classical case is that our variations need to preserve the filtrations, and our processes will be adapted.  If  $(\Omega,\mathcal{A}, (\mathcal{A}_t)_{t\in[0,1]},\mathcal{P})$ is a complete stochastic basis (see above) we  denote by $L^2_a(\mathcal{P},H)$  (similarly for the other $L^p(\mathcal{P}, H)$ and  $L^p(\mathcal{P}, H_{0,0})$ spaces) the subspace of  $u\in L^2(\mathcal{P},H)$ such that $(t,\omega) \to   u_t(\omega)\in \mathbb{R}^d$ is $(\mathcal{A}_t)-$adapted for any (and then all) continuous processes whose $\mathcal{P}-$equivalence class is $u$.  For $u \in L^0_a(\mathcal{P}, W_{abs})$ we can always find  $\int_0^. \dot{u}_s ds$  in the equivalence class of $u$ so that $(t,\omega) \to \dot{u}_t$ is $(\mathcal{A}_t)-$predicable: we choose such modifications of the derivative unless  expressively stated. 
In particular, for $\eta\in \mathcal{P}_W$, $(\mathcal{F}_t^\eta)_{t\in[0,1]}$ will denote the $\eta-$usual augmentation of the filtration generated by the evaluation process $(W_t)$, with the convention $\mathcal{F}_1^\eta= \mathcal{B}(W)^\eta$. The space $L^2_a(\eta, H)$ (similarly for the other $L^p_a$ spaces) will denote the set $L^2_a(\mathcal{P},H)$ for $(\Omega, \mathcal{A}, \mathcal{P})= (W,\mathcal{B}(W)^\eta, \eta)$ and the filtration $(\mathcal{F}_t^\eta)$. In the whole paper $\lambda$ will denote the Lebesgue measure on $[0,1]$. Finally, for convenience of notations, if $u:= \int_0^. \dot{u}_s ds\in L^0_a(\mathcal{P},W_{abs})$, and $\tau$ is an $(\mathcal{A}_t)-$stopping time we note $$\pi_\tau u :=\int_0^{.\wedge \tau}\dot{u}_s ds$$  the process stopped by $\tau$.

\subsection{Martingales by duality}

The variational characterization of martingales is a result of stochastic control (see \cite{EMERY} and the references therein) which relies on duality. Since it will play a central role  in this paper we provide here a precise statement of this result. 

Let $(\Omega, \mathcal{A}, (\mathcal{A}_t)_{t\in[0,1]}, \mathcal{P})$ be a complete stochastic basis. The mapping $$r :  \beta \in  L^2(\mathcal{P},\mathbb{R}^d)   \to \int_0^. E_\nu\left[\beta \middle| \mathcal{F}_s^\nu\right] ds  \in L^2_a(\mathcal{P},H) $$  defines a linear operator which is continuous by Jensen's inequality. Its adjoint is given by the operator $$q : u \in L^2_a(\mathcal{P},H) \to u_1:=\int_0^1 \dot{u}_s ds \in L^2(\mathcal{P},\mathbb{R}^d) $$  which is also linear and continuous. Indeed from the definitions we obtain directly \begin{equation}  E_{\mathcal{P}}[<q(u), \beta>_{\mathbb{R}^d}] = E_{\mathcal{P}}[<u, r(\beta)>_{H}] \end{equation} for any $\beta\in L^2({\mathcal{P}},\mathbb{R}^d)$ and $u\in L^2_a({\mathcal{P}},H)$.  By a classical result of functional analysis (see \cite{RR} Chapter  VI Lemma 6 for instance)  the orthogonal  of the kernel of $q$ (i.e. $q^{-1}(\{0_{L^2({\mathcal{P}},\mathbb{R}^d)}\})$) in the Hilbert space $L^2_a({\mathcal{P}},H)$ coincides with the closure of the range of $r$ in $L^2_a({\mathcal{P}},H)$. As a matter of fact, by a stopping argument,  it is straightforward to see that this latter space is  the space of maps $u\in L^2_a({\mathcal{P}},H)$ with a martingale derivative. A precise statement of this result is the following orthogonal decomposition of $L^2_a({\mathcal{P}},H)$ which immediatly yields the variational characterization of the martingale:

\begin{proposition}
\label{emery}
For any complete stochastic basis $(\Omega, \mathcal{A}, (\mathcal{A}_t)_{t\in[0,1]},\mathcal{P})$ we have \begin{equation} \label{emeryorth}L^2_a(\mathcal{P},H) = \mathcal{M}_a(\mathcal{P}, H) \oplus_\perp L^2_a(\mathcal{P}, H_{0,0}) \end{equation} where $\mathcal{M}_a(\mathcal{P}, H)$ is the set of  $u\in L^2_a({\mathcal{P}},H)$ for which there exists a  \textit{c\`{a}dl\`{a}g} $(\mathcal{A}_t)-$martingale $(M_t)_{t\in[0,1)}$ such that ${\mathcal{P}}-$a.s.
$$u= \int_0^. M_s ds$$ and where \begin{equation} \label{knu} L^2_a(\mathcal{P}, H_{0,0})  := \{ h \in L^2_a({\mathcal{P}},H) :  {\mathcal{P}}-a.s.~ h_1= 0 \}\end{equation}
In particular, for  $ \alpha \in L^2({\mathcal{P}}, \mathbb{R}^d)$, if $$\mathcal{C}(\alpha):=\left\{u\in L^2_a({\mathcal{P}},H), \mathcal{P}-a.s. ~u_1=\alpha\right \}$$ and $$I : \alpha \in L^2({\mathcal{P}}, \mathbb{R}^d) \to I(\alpha) \in \mathbb{R} \cup\{ \infty\}$$ is defined by $$I(\alpha):= \inf\left(\left\{E_{{\mathcal{P}}}\left[|u|_H^2\right] : u\in \mathcal{C}(\alpha) \right\}\right) ,$$ for any $\alpha\in D_I:=\{ \alpha\in L^2({\mathcal{P}}, \mathbb{R}^d)  : I(\alpha)<\infty \}$ the infimum is attained by a unique element $u^\star(\alpha)\in \mathcal{C}(\alpha)$, which is the orthogonal projection of any (and then of all)  element(s) of $\mathcal{C}(\alpha)$ on $\mathcal{M}_a(\mathcal{P}, H)$. Conversely if a $u\in \mathcal{C}(\alpha)$ is an element of $\mathcal{M}_a(\mathcal{P}, H)$ it attains the infinimum of $I(\alpha)$. 
\end{proposition}  
\begin{remarkk}
\label{emeryrem}
For convenience of notations we considered $\mathbb{R}^d-$valued processes in the proofs of  Proposition~\ref{emery} and Proposition~\ref{emerydual} but the result also holds for processes with values in any separable Hilbert space, as it is well known. Moreover by taking some trivial probability space   one obtains as a particular case that the orthogonal to $H_{0,0}:=\{h\in H : h_1 =0\}$ in $H$ is the set of  $h\in H$ such that there exists a $c_h\in \mathbb{R}^d$ with a.s. for all $s\in[0,1]$ $\dot{h}_s= c_h$.

\end{remarkk}


%

The following result is dual to Proposition~\ref{emery} :

\begin{proposition}
\label{emerydual}
Let $(\Omega,\mathcal{A}, (\mathcal{A}_t)_{t\in[0,1]}, \mathcal{P})$ be a complete stochastic basis and let $L^2(\mathcal{P}, \mathbb{R}^d ; \mathcal{A}_{1-})$ be the set of the $\alpha \in L^2(\mathcal{P}, \mathbb{R}^d)$  such that $\alpha$ is $\mathcal{A}_{1-}$ measurable. Then the set of  $\alpha\in L^2_a(\mathcal{P},\mathbb{R}^d)$ that can be attained by an adapted shift (i.e. such that there exists a  $u\in L^2_a({\mathcal{P}}, H)$ with ${\mathcal{P}}-a.s.$ $u_1=\alpha$)  is a dense subspace of  $L^2(\mathcal{P}, \mathbb{R}^d ;  \mathcal{A}_{1-})$ for the $L^2({\mathcal{P}}, \mathbb{R}^d)$ topology.
\end{proposition}
\nproof First note that the set of  $\alpha\in L^2({\mathcal{P}},\mathbb{R}^d)$ that can be attained by an adapted shift coincides with the range $q(L^2_a({\mathcal{P}}, H))$  of $q$. Hence, if we denote by $\overline{q(L^2_a({\mathcal{P}}, \mathbb{R}^d))}$ the closure of  $q(L^2_a({\mathcal{P}}, H))$, we have to prove that $ \overline{q(L^2_a({\mathcal{P}}, \mathbb{R}^d))} = L^2({\mathcal{P}},\mathbb{R}^d; \mathcal{A}_{1-})$.   By continuity  $q(L^2_a({\mathcal{P}}, \mathbb{R}^d) \subset L^2(\mathcal{P}, \mathbb{R}^d ;  \mathcal{A}_{1-})$, and since this latter space is closed we obtain  \begin{equation}\label{incl1}  \overline{q(L^2_a({\mathcal{P}}, \mathbb{R}^d))} \subset L^2({\mathcal{P}},\mathbb{R}^d; \mathcal{A}_{1-})   \end{equation} We now prove the converse inclusion.  By duality, the closure of $q(L^2_a({\mathcal{P}}, H))$  is the orthogonal in $L^2_a({\mathcal{P}},H)$ to the kernel $r^{-1}(\{0\})$ of $r$, which is given by $$ r^{-1}(\{0\}) := \left\{ \alpha\in L^2({\mathcal{P}},\mathbb{R}^d): {\mathcal{P}}-a.s.\ \int_0^. E_{\mathcal{P}}[\alpha| \mathcal{A}_s] ds = 0 \right\}. $$ By considering a right continuous modification of  $(E_{\mathcal{P}}[\alpha| \mathcal{A}_t])_{t\in[0,1]}$,  the martingale convergence theorem yields \begin{equation} \label{11M} r^{-1}(\{0\})= \{ \mathcal{P}-a.s. \ E_{\mathcal{P}}[\alpha| \mathcal{A}_{1-}] = 0 \}.\end{equation} Let $X \in L^2({\mathcal{P}},\mathbb{R}^d;\mathcal{A}_{1-})$ and $\alpha\in r^{-1}(\{0\})$. Then, by definition, $$E_{\mathcal{P}}\left[<X,\alpha>_{\mathbb{R}^d}\right] = E_{\mathcal{P}}\left[<X,E_{\mathcal{P}}\left[ \alpha \middle|\mathcal{A}_{1-}\right]>_{\mathbb{R}^d}\right]=0.$$  Hence $$ L^2({\mathcal{P}},\mathbb{R}^d; \mathcal{A}_{1-}) \subset r^{-1}(\{0\})^\perp =   \overline{q(L^2_a({\mathcal{P}}, \mathbb{R}^d))}$$ Together with~(\ref{incl1}) we obtain the desired result. \nqed

\subsection{Transformations of measure preserving the filtration}
\label{tpf}

In this section we introduce isomorphisms of a filtered probability space, which are usually used to perform transformations of measure preserving the filtrations, in particular in Malliavin calculus. Here we will handle morphisms of probability spaces (see above). Indeed the results we use only provide  existence of equivalence classes of mappings measurable with respect to completed sigma fields. Recall that $M_{\mathcal{P}}((\Omega, \mathcal{A}),(W,\mathcal{B}(W))$ denotes the set of $\mathcal{P}-$equivalence class of $\mathcal{A}/ \mathcal{B}(W)-$measurable mappings $f : \Omega \to W$. To avoid heavy notations, whenever we  handle a property which does not depend on the element in the equivalence class, we implicitly denote with the same letter  $U\in M_{\mathcal{P}}((\Omega, \mathcal{A}),(W,\mathcal{B}(W))$  and  a $\mathcal{A}/ \mathcal{B}(W)-$measurable mapping in this class. However within this whole subsection  we will make the difference, in order to avoid any ambiguity on the notations. The main properties related to transformations of measure preserving the filtrations concern their inverse images  and pullbacks.   If $\mathcal{G}$ is a sigma-field and $U$ is a $\mathcal{P}-$equivalence class of $\mathcal{A}/ \mathcal{B}(W)-$measurable mappings (i.e. $U\in M_{\mathcal{P}}((\Omega, \mathcal{A}),(W,\mathcal{B}(W))$), we denote by $U^{-1}(\mathcal{G})$ the $\mathcal{P}-$completion of $f^{-1}(\mathcal{G})$ for any (and then all) $\mathcal{A}/ \mathcal{B}(W)-$measurable $f : \Omega \to W$ such that $\mathcal{P}-a.s.$ $U=f$ (i.e. $U$ is the $\mathcal{P}-$equivalence class of $f$, see above) and we call it the inverse image of $\mathcal{G}$ by $U$. This name is justified by its behaviour by pullback which we now recall.

  Given  $\eta\in \mathcal{P}_W$, and $U$  (resp. $X$) a $\eta-$equivalence class of $\mathcal{B}(W)^\eta / \mathcal{B}(W)-$measurable mappings (resp. a $\mathcal{P}-$equivalence class of $\mathcal{A}/ \mathcal{B}(W)-$measurable mappings),  under the assumption that $X_\star\mathcal{P}$ is absolutely continuous w.r.t. $\eta$ (i.e. $X_\star\mathcal{P}<<\eta$)  we have $\mathcal{P}-a.s.$ $$f_U\circ g_X = f_{\widetilde{U}}\circ g_{\widetilde{X}}$$  for any measurable $f_U, f_{\widetilde{U}} : W\to W$ (resp. $g_X,g_{\widetilde{X}} : \Omega \to W$) in the $\eta-$equivalence class $U$ (resp.  in the $\mathcal{P}-$equivalence class $X$), where $f_U\circ g_X : \omega \in \Omega \to f_U(g_X(\omega))\in W$ (similarly for $ f_{\widetilde{U}}\circ g_{\widetilde{X}}$). We denote by $U\circ X$ the $\mathcal{P}-$equivalence class of the $\mathcal{A}/ \mathcal{B}(W)$ measurable mapping $f_U\circ g_X$ for any (and then all) such $f_U$ and $g_X$. Then, for all sigma field $\mathcal{G}$  of $W$ \begin{equation} \label{ucircX} (U\circ X)^{-1}(\mathcal{G}) = X^{-1}(U^{-1}(\mathcal{G})). \end{equation}  
  
    This is related to adapted processes in the following way. Denote by $(\mathcal{B}_t^0(W))_{t\in [0,1]}$ the filtration generated by the evaluation process on $W$ i.e. for all $t\in [0,1]$ $$\mathcal{B}_t^0(W):= \sigma(W_s, s\leq t),$$  Since we shall deal with
 progressively measurable processes and \textit{c\`{a}dl\`{a}g} modifications of martingales, for $\eta\in \mathcal{P}_W$ we will consider its usual augmentation $(\mathcal{F}_t^\eta)_{t\in [0,1]}$ (under $\eta$). We recall that  $$\mathcal{F}_t^\eta := \mathcal{B}_{t+}^0(W)^\eta$$ for $t\in [0,1]$. Here we adopt the conventions that at $t=1$ the usual augmentation is just the completion and that $\mathcal{B}_{1+}^0(W):= \mathcal{B}(W)$. Similarly, for  $U\in M_{\mathcal{P}}((\Omega, \mathcal{A}),(W,\mathcal{B}(W))$, we will need to consider the following filtration generated by $U$. To any $f: \Omega \to W$ which is $\mathcal{A}/ \mathcal{B}(W)-$measurable  for all $t\in [0,1]$  denote  $$\mathcal{G}_t^f := \sigma( f_s, s\leq t)$$ where $(f_t)$ is the measurable process associated to $f$ by $$(f_t) : (t,\omega) \in [0,1] \times W \to f_t(\omega):=W_t(f(\omega)) \in \mathbb{R}^d .$$ Note that by definition we also have for all $t\in[0,1]$, $\mathcal{G}_t^f = f^{-1}(\mathcal{B}^0_t(W))$, and that it is elementary to check that $\mathcal{G}_{t+}^f= f^{-1}(\mathcal{B}^0_{t+}(W))$. Then, if  $(\Omega, \mathcal{A}, (\mathcal{A}_t)_{t\in[0,1]}, \mathcal{P})$ is a complete stochastic basis we  say  that $U$ is $(\mathcal{A}_t)-$adapted if and only if any (and then all) $\mathcal{A}/ \mathcal{B}(W)-$measurable $f: \Omega \to W$ such that $\mathcal{P}-a.s.$ $U= f$ is $(\mathcal{A}_t)-$adapted i.e. for all $t\in [0,1]$, $\mathcal{G}_t^f \subset \mathcal{A}_t$.  We define the filtration generated by $U$, which we note $(\mathcal{G}_t^U)$,  to be the usual augmentation with respect to $\mathcal{P}$ of the filtration $(\mathcal{G}_t^f)$ for any (and then all) $\mathcal{A}/ \mathcal{B}(W)-$measurable $f$ such that $\mathcal{P}-a.s.$ $U= f$. In particular for all $t\in [0,1]$  it is elementary to check that with these definitions \begin{equation} \label{gubet} \mathcal{G}_t^U = U^{-1}(\mathcal{B}^0_{t+}(W)) = U^{-1}(\mathcal{F}_t^{U_\star \mathcal{P}})\end{equation} and that, due to our hypothesis on $(\mathcal{A}_t)$, $U$ is $(\mathcal{A}_t)$-adapted if and only if for all $t\in [0,1]$ \begin{equation} \label{gualp} \mathcal{G}_t^U \subset \mathcal{A}_t .\end{equation} Thus, by~(\ref{ucircX}) if $U$ is $(\mathcal{A}_t)$-adapted and $X$ is $(\mathcal{F}_t^\eta)$ adapted $U\circ X$ is also $(\mathcal{A}_t)$-adapted. Conversely, an easy criterion for the existence of an adapted pullback is the following Proposition. We emphasize  that it only yields the existence of a measurable function which is measurable w.r.t. the completed space with equality up to  negligible sets.

 \begin{proposition}
\label{propadapt}
Assume that $Y, X$ are two $\mathcal{P}-$equivalence classes of $\mathcal{A}/ \mathcal{B}(W)-$ measurable mappings (i.e. two elements of $M_{\mathcal{P}}((\Omega,\mathcal{A}), (W,\mathcal{B}(W))$. Then the following assertions are equivalent 
\begin{enumerate}[(i)]
\item $Y$ is adapted to the filtration generated by $X$ i.e. for all $t\in [0,1]$ $$\mathcal{G}_t^Y \subset \mathcal{G}_t^X$$ where $(\mathcal{G}_t^X)$ (resp. $(\mathcal{G}_t^Y)$) is the $\mathcal{P}-$usual augmentation of the filtration generated by any $\mathcal{A}/\mathcal{B}(W)$ mesurable  $f_X : \Omega \to W$ (resp. $g_Y: \Omega \to W$) whose $\mathcal{P}-$equivalence class is $X$ (resp. $Y$).
\item There exists a  $F\in M_{X_\star \mathcal{P}}((W,\mathcal{B}(W)^{X_\star \mathcal{P}}),(W,\mathcal{B}(W))$ which is  $(\mathcal{F}_t^{X_\star \mathcal{P}})-$ adapted such that $\mathcal{P}-a.s.$  $$Y= F \circ X $$ where $F\circ X $ denotes the pullback defined above, and $(\mathcal{F}_t^{X_\star \mathcal{P}})-$is the $X_\star \mathcal{P}$-usual augmentation of the filtration generated by the evaluation process. In particular $$F_\star(X_\star \mathcal{P})= Y_\star\mathcal{P}$$
\end{enumerate}
Moreover the two following assertions are equivalent :

\begin{enumerate}[(1)]
\item $X$ and $Y$  generate the same filtrations i.e. for all $t\in [0,1]$ $$\mathcal{G}_t^Y = \mathcal{G}_t^X$$ 
\item There exists a  $F\in M_{X_\star \mathcal{P}}((W,\mathcal{B}(W)^{X_\star \mathcal{P}}),(W,\mathcal{B}(W))$ which is  $(\mathcal{F}_t^{X_\star \mathcal{P}})-$ adapted
and a $G\in M_{Y_\star \mathcal{P}}((W,\mathcal{B}(W)^{Y_\star \mathcal{P}}),(W,\mathcal{B}(W))$ which is $(\mathcal{F}_t^{Y_\star \mathcal{P}})-$ adapted such that
 $\mathcal{P}-a.s.$  $$Y= F \circ X $$ and $$X = G\circ Y$$  Moreover $X_\star \mathcal{P}-a.s.$ $$G\circ F = I_W$$ and $Y_\star \mathcal{P}-a.s.$ $$F\circ G = I_W$$
\end{enumerate}

\end{proposition}
\nproof 
 Similar to the proof of the Yamada-Watanabe criterion (see \cite{I-W}).
 \nqed

 The isomorphisms of filtered probability spaces play a key role in Malliavin's work. We now state their definition.  First note that whenever $\mathcal{A}$ is complete, $f : \Omega \to W$ is $\mathcal{A}/ \mathcal{B}(W)-$measurable if and only if it is $\mathcal{A}/ \mathcal{B}(W)^{f_\star \mathcal{P}}$ measurable. This ensures that the pullbacks below are well defined. Let $\eta, \nu \in \mathcal{P}_W$. We say that $U \in M_{\mathcal{P}}((W, \mathcal{B}(W)^\eta),(W,\mathcal{B}(W))$ with $U_\star \eta = \nu$  is an isomorphism of filtered probability spaces on $(W, \mathcal{F}_.^\eta, \eta)$ to $(W, \mathcal{F}_.^\nu, \nu)$  if  $U$ is $(\mathcal{F}_t^\eta)-$adapted and if there exists  $\widetilde{U}\in M_{\nu}((W,\mathcal{B}(W)^\nu), (W, \mathcal{B}(W))$ which is $(\mathcal{F}_t^\nu)-$adapted and such that $\eta-a.s.$ $$\widetilde{U}\circ U =I_W$$ and $\nu-a.s.$ $$U\circ \widetilde{U} = I_W$$ where $I_W : \omega\in W \to \omega \in W$. In this case $\widetilde{U}$ is unique and we call it the inverse of $U$. Note that by~(\ref{ucircX}) and~(\ref{gubet}) we have,  $$\mathcal{G}_t^U = \mathcal{F}_t^\eta$$ and $$\mathcal{G}_t^{\widetilde{U}} = \mathcal{F}_t^\nu$$ for all $t\in [0,1]$. We emphasize that, with this definition, the invertibility may fail on some negligible set. Explicitly for $f_U$ (resp. $f_{\widetilde{U}}$) $\mathcal{B}(W)^\eta/ \mathcal{B}(W)$ (resp. $\mathcal{B}(W)^\nu / \mathcal{B}(W)$) measurable whose equivalence class is $U$ (resp. $\widetilde{U}$) the equality $f_U\circ f_{\widetilde{U}}=I_W$ (resp. $f_{\widetilde{U}}\circ f_{U}=I_W$) is merely assumed to hold $\eta-a.s.$ (resp. $\nu-a.s.$). In particular this definition doesn't claim the existence of a Borel measurable bijection (invertible everywhere) in the equivalence class of $U$.   Note also that any such $U$ induces an obvious isometric identification of the $L^p$ spaces of $\eta$ and $U_\star \eta $ (and of the $L^p_a(\eta, H)$ and $L^p_a(U_\star\eta ,H))$ which is sometimes used as an alternative definition (just consider $k\in L^p_a(U_\star \eta,H) \to k\circ U \in L^p_a(\eta,H)$ for instance). Useful characterizations  to handle isomorphisms of filtered probability spaces are  provided by the following proposition:

\begin{proposition}
\label{propisom} 
Let $\eta \in \mathcal{P}_W$ and $U\in M_{\mathcal{P}}((W, \mathcal{B}(W)^\eta),(W,\mathcal{B}(W))$. Then the following  are equivalent
\begin{enumerate}[(i)]
\item $U$ is an isomorphism of filtered probability spaces on  $(W,\mathcal{F}_.^\eta,\eta)$ to $(W,\mathcal{F}_.^{U_\star \eta}, U_\star\eta)$ 
\item  For all  $t\in [0,1]$ \begin{equation} \label{filtreg} \mathcal{G}_t^U = \mathcal{F}_t^\eta\end{equation}
\item~(\ref{filtreg}) holds for all $t\in [0,1)$ and, for all $X, Y\in M_{\mathcal{P}}((\Omega, \mathcal{A}),(W,\mathcal{B}(W))$ 
defined on the same complete space $(\Omega, \mathcal{A}, \mathcal{P})$  such that $$X_\star \mathcal{P} = Y_\star \mathcal{P} = \eta ,$$ we have $$\mathcal
{P}-a.s. \ U(X) = U(Y) \implies \mathcal{P}-a.s.  ~ X= Y $$
\item~(\ref{filtreg}) holds for all $t\in [0,1)$  and, for every  complete probability space $(\Omega,\mathcal{A},\mathcal{P})$ and  for all $Y \in M_P((\Omega, \mathcal{A}), (W, \mathcal{B}(W))$ such that $Y_{\star}\mathcal{P} = U_\star \eta$, there exists a $X \in M_P((\Omega, \mathcal{A}), (W, \mathcal{B}(W))$ with $X_\star \mathcal{P} = \eta$ such that $\mathcal{P}-a.s.$ $$Y = U\circ X$$
\end{enumerate}
Moreover, in the case where one of the above assumptions is satisfied, $X$ in $(iv)$ is unique.
\end{proposition}
\nproof
 Similar to the proof of the Yamada-Watanabe criterion (see \cite{I-W}). 
 \nqed
\begin{remarkk}
In practice $(ii)$ is useless to obtain $(i)$; since one may use $(iii)$ to prove that $\mathcal{G}_1^U = \mathcal{F}_1^\eta = \mathcal{B}(W)^\eta$, $(iii)$ is the most efficient criterion to obtain $(i)$.
\end{remarkk}

In the sequel we will need to control the initial behaviour of isomorphisms, namely we will need them to preserve the initial information. For this reason we set the following definition.
\begin{definition}
\label{iodef}
For $\eta\in \mathcal{P}_W$, we denote by $\mathcal{I}_f^0(\eta)$ the set of  isomorphisms of filtered probability spaces  $U\in M_{\eta}((W, \mathcal{B}(W)^\eta),(W,\mathcal{B}(W))$ on  $(W,\mathcal{F}_.^\eta,\eta)$ to $(W,\mathcal{F}_.^{U_\star \eta}, U_\star\eta)$ which further satisfy  $$\sigma(W_0)^\eta = \sigma(U_0)^\eta$$ 
\end{definition}
\begin{remarkk}
Note that by Dynkin's lemma an isomorphism of probability spaces  on $(W,\mathcal{F}_.^\eta,\eta)$ to $(W,\mathcal{F}_.^{U_\star \eta}, U_\star\eta)$ is an element of $\mathcal{I}_f^{0}(\eta)$ if and only if  there exists a $\mathcal{B}(\mathbb{R}^d)^{{U_0}_\star \eta} / \mathcal{B}(\mathbb{R}^d)-$ measurable function $f : \mathbb{R}^d\to \mathbb{R}^d$  and a $\mathcal{B}(\mathbb{R}^d)^{{W_0}_\star\eta}/ \mathcal{B}(\mathbb{R}^d)-$ measurable function $g : \mathbb{R}^d\to \mathbb{R}^d$ satisfying ${U_0}_\star \eta-a.s.$ (resp. ${W_0}_\star\eta-a.s. $) $g\circ f =I_{\mathbb{R}^d}$ (resp. $f\circ g =I_{\mathbb{R}^d}$) such that $\eta-a.s.$ $ U_0= g(W_0)$ and $U_\star \eta-a.s.$ $\widetilde{U}_0= f(W_0)$ where $\widetilde{U}$ denotes the inverse of $U$. In particular for $\nu \in \mathcal{P}_W$, if $U$ is an isomorphism of probability spaces on $(W, \mathcal{F}_.^\eta, \eta)$ to $(W, \mathcal{F}_.^\nu, \nu)$ then $U\in \mathcal{I}_f^0(\eta)$ if and only if $\widetilde{U}\in \mathcal{I}_f^0(\nu)$. Using this it is straightforward to check that if $U\in \mathcal{I}_f^0(\eta)$ and $T\in \mathcal{I}_f^0(U_\star\eta)$ then $T\circ U \in \mathcal{I}_f^0(\eta)$.
\end{remarkk}

\subsection{Some spaces of laws of continuous semi-martingales}

\subsubsection{The space $\mathbb{S}$}

Within our stochastic extensions, the space $\mathbb{S}$ will play a role analogous to the path space $\Omega_{[0,1]}^2$ of the $\mathbb{R}^d$-valued $C^2-$functions on $[0,1]$ in the classical calculus of variations. In the whole paper $\mathbb{S}$ will denote the space of the Borel probabilities $\nu \in \mathcal{P}_W$ for which there exist 

\begin{enumerate}[(i)]
\item  A continuous $(\mathcal{F}_t^\nu)-$ local martingale $(\widetilde{M}_t^\nu)$ and a $(\mathcal{F}_t^\nu)-$ predicable process $(\widetilde{v}_t^\nu)$ defined on the space $(W,\mathcal{B}(W)^\nu,\nu)$ such that $\nu-a.s.$ for any $t\in [0,1]$  \begin{equation} W_t = W_0 + \widetilde{M}_t^\nu + \int_0^t \widetilde{v}_s^\nu ds \end{equation} 
\item Two $\mathcal{M}_d(\mathbb{R})$ valued $(\mathcal{F}_t^\nu)$ predicable processes $(\widetilde{\alpha}^\nu_t)$ and $(\widetilde{\sigma}^\nu_t)$ related by $$\widetilde{\alpha}^\nu_t := \widetilde{\sigma}^\nu_t {(\widetilde{\sigma}^\nu_t)}^{\dagger}$$ such that $$\int_0^1 |(\widetilde{\alpha}_t^\nu)^{ij}| dt < \infty$$ $$\int_0^1 |(\widetilde{\sigma}_t^\nu)^{ij}|^2 dt <\infty$$ for all $i,j\in [1,d]$ and  $\nu-a.s.$  $$<(\widetilde{M}^{\nu})^i, (\widetilde{M}^{\nu})^j> = \int_0^. (\widetilde{\alpha}^\nu_s)^{ij} ds$$ where $\mathcal{M}_d(\mathbb{R})$  denotes the set of  $d\times d$ matrices endowed with its usual topology, and where $<.,.>$ denotes the predicable quadratic co-variation process.
\end{enumerate}

Note that for $\nu \in \mathbb{S}$ the continuous local martingale $(\widetilde{M}_t^\nu)$ and the finite variation term $(\int_0^. \widetilde{v}_t^\nu dt)$ are unique up to a $\nu$-evanescent set. Hence to $\nu$ we associate canonically its martingale part (resp. its finite variation part) which is defined to be $M^\nu$ (resp. $b^\nu$), the element of $L^0_a(\nu, W)$ (resp.  $L^0_a(\nu, W_{abs})$)     such that $\nu-a.s.$ for all $t$, $W_t\circ M^\nu = \widetilde{M}_t^\nu$ (and $\nu-a.s.$ $b^\nu=\int_0^. \widetilde{v}_t^\nu dt$). On the other hand we note $(v^\nu_t)$ (resp. $(\alpha^\nu_t)$) the equivalence classes  of $(\mathcal{F}_t^\nu)-$optional processes which are  $\lambda\otimes \nu$ equal to $(\widetilde{v}_t^\nu)$ (resp. to $(\widetilde{\alpha}_t^\nu)$). Actually we can always chose a $(\mathcal{F}_t^\nu)-$predicable process in the equivalence class of  $(v_t^\nu)$ (resp. of $(\alpha_t^\nu)$) and we will do this, unless it is explicitly stated that we take it optional, or more precisely (when such a modification exists) that we take it right continuous and $(\mathcal{F}_t^\nu)-$adapted.

\begin{definition}
\label{defvelocity}
For $\nu\in \mathbb{S}$ we note $M^\nu$ (resp. $b^\nu$), both elements of $M_{\nu}((W, \mathcal{B}(W)^\nu), (W, \mathcal{B}(W))$, the martingale part (resp. the finite variation part) of $\nu$ and we call $(\int_0^. v_t^\nu dt)$ (resp. $(\int_0^. \alpha_t^\nu dt)$), such that $\nu-a.s.$ $$<M^\nu>= \int_0^. \alpha_t^\nu dt$$ and $\nu-a.s.$ $$b^\nu = \int_0^. v_t^\nu dt$$ the characteristics of $\nu$.
\end{definition}

\section{Variation processes}
\label{section2}
\subsection{Vector space of variation processes}

\begin{definition}
\label{defframe}
For $\nu \in \mathcal{P}_W$ we denote by $V_\nu$  the set of  $h\in L^2_a(\nu,H)$ such that for any $U\in \mathcal{I}_f^0(\nu)$ (see Definition~\ref{iodef}) we also have $U^{h}\in \mathcal{I}_f^0(\nu)$ where  \begin{equation}\label{thetah} U^{h}:= U +  h \end{equation} Moreover $V_\nu^\infty$ (resp. $V_\nu^0$, resp. $V_\nu^{0,\infty}$) will denote the linear subspace of $V_\nu$ defined by $$V_\nu^\infty :=  V_\nu \cap L^\infty(\nu,W) $$
$$V_\nu^0:= V_\nu \cap L^2(\nu, H_{0,0})$$
$$V_\nu^{0,\infty}:= V_\nu^\infty  \cap V_\nu^0$$
 We  say that  $h\in V_\nu$ is the variation process of the curve $({\tau_{\epsilon h}}_\star \nu)_{\epsilon\in \mathbb{R}}\subset \mathcal{P}_W$ at $\nu$ where for any $k\in L^2_a(\nu,H)$, $\tau_k$ denotes $$\tau_k := I_W + k$$   \end{definition}

\begin{proposition}
 For any $\nu \in \mathcal{P}_W$, $V_\nu$ is a vector subspace of $L^2_a(\nu,H)$. 
\end{proposition}
\nproof
Consider $\nu \in \mathcal{P}_W$, an element $h\in V_\nu$, $\epsilon \in \mathbb{R}$. We first prove that $\epsilon h \in V_\nu$.  Let $U\in \mathcal
{I}_f^0(\eta)$. For  $\epsilon\in \mathbb{R} / \{0\}$ set $$U^{\epsilon h} = U + \epsilon h$$ We have to prove that $U^{\epsilon h}\in \mathcal{I}_f^0
(\eta)$.   Note that \begin{equation} \label{uepsh} U^{\epsilon h} = \epsilon (U^\epsilon + h) \end{equation}  where $$U^\epsilon := \frac{1}{\epsilon} U$$ Since $U\in \mathcal{I}_f^0(\nu)$ by Proposition~\ref{propisom} we also have $U^\epsilon\in \mathcal{I}_f^0(\nu)$. Hence by the definition of $V_\nu$, 
$U^\epsilon + h \in \mathcal{I}_f^0(\nu)$. Similarly, by~(\ref{uepsh}) we have $U^{\epsilon h}\in \mathcal{I}_f^0(\nu)$.  Therefore, if $h\in V_\nu$ 
and $\epsilon \in \mathbb{R}$, $\epsilon h \in V_\nu$. We take $h,k\in V_\nu$ and we want to derive $h+k\in V_\nu$. For $U\in \mathcal{I}
_f^0(\nu)$, we need to prove that $$U^{k+h}:= U +  h + k \in \mathcal{I}_f^0(\nu).$$ Since $k\in V_\nu$,  $U^{k}:= U +  k \in \mathcal{I}_f^0(\nu)$. 
On the other hand since $h\in V_\nu$ and  $U^{h+k}= U^k +  h$, the definition of $V_\nu$ yields $U^{h+k}\in \mathcal{I}_f^0(\nu)$. Therefore  $V_\nu
$ is a vector space.  
\nqed

\begin{remarkk}
First note that for $\nu\in \mathcal{P}_W$, by Proposition~\ref{propisom}, it is straightforward to check that  $$H \hookrightarrow V_\nu^\infty \subset V_\nu  \subset L^2_a(\nu,H)$$ and $$H_{0,0} \hookrightarrow  V_\nu^{\infty,0} \subset V_\nu^{0} \subset L^2_a(\nu,H_{0,0}).$$ In the case where  $\nu=\delta^{Dirac}_h$ (Dirac measure concentrated on $h\in H$) we have $L^2_a(\nu,H)\simeq H$ and $L^2_a(\nu,H_{0,0})\simeq H_{0,0} $ so that all the inclusions become equalities. On the other hand if we take $\nu\in \mathbb{S}$ to be the law of the solution to Tsirelson's equation (see \cite{Tsi} or \cite{I-W})  which we note $$dX_t= dB_t +v_t(X) dt $$  then the $\nu$-equivalence class of $I_W$ is an isomorphism of probability spaces on $(W,\mathcal{F}^\nu_.,\nu)$ and $\int_0^. v_s ds\in L^2_a(\nu,H)$ but $I_W -\int_0^. v_s  ds$ is not an isomorphism of filtered probability spaces (see \cite{RL-locinv}, \cite{ASU-3},\cite{RLASU}). By localization  one can build examples of probabilities in $\mathbb{S}$ where $V_\nu^{0}$ is a proper linear subspace of $L^2_a(\nu,H_{0,0})$. However we shall see that fortunately for any $\nu\in \mathcal{P}_W$ these injections  (except those of $H$ and $H_{0,0}$) are always dense in the topology of $L^2(\nu,H)$ for any $\nu \in \mathbb{S}$.
\end{remarkk}

The following Proposition shows that the variation processes are invariant by isomorphisms. Since we will not use it in the sequel, it can be skipped in a first reading. 

\begin{proposition}
\label{spapropisom}
 For any $\nu \in \mathcal{P}_W$ and any $U \in \mathcal{I}_f^0(\nu)$ we have  $$V_{U_\star \nu} \simeq_{j_U} V_\nu$$
More precisely the mapping $j_U$ defined by $j_U : h\in V_{U_\star \nu} \to h\circ U \in V_\nu$ is a bijection (and an isometry) of $V_{U_\star \nu}$ onto $V_\nu$  whose inverse is given by 
 $j_{\widetilde{U}} : h\in V_\nu \to h\circ\widetilde{U} \in  V_{U_\star \nu}$  where $\widetilde{U}$ is the inverse of $U$.
\end{proposition}
\nproof
Consider $\nu \in \mathcal{P}_W$  and $U\in \mathcal{I}_f^0(\nu)$, whose inverse is denoted by $\widetilde{U}\in \mathcal{I}_f^0(U_\star\nu)$. By  symmetry, to prove the result it is sufficient to swhow that for $h\in V_{U_\star \nu}$ we have  $j_U(h) \in V_{\nu}$. Hence we consider $T \in  \mathcal{I}_f^0(\nu)$, whose inverse is denoted by $\widetilde{T} \in \mathcal{I}_f^0(T_\star\nu)$ and we prove that  $$T^{j_U(h)}:= T+ j_U(h) $$ is an element of $\mathcal{I}_f^0(\nu)$. To see this, note that  $T\circ \widetilde{U}\in \mathcal{I}_f^0(U_\star\nu)$ (with inverse $U\circ \widetilde{T} \in \mathcal{I}_f^0(T_\star\nu)$). Since $h\in V_{U_\star \nu}$, we have $$T^h:=T\circ \widetilde{U} + h \in \mathcal{I}_f^0(U_\star\nu)$$ We denote by $S^h\in \mathcal{I}_f^0((T^h\circ U)_\star\nu)$ the inverse of $T^h$. Finally, by definition, $\nu-a.s.$ $$T^{j_U(h)}= T^h\circ U $$ so that $T^{j_U(h)}\in \mathcal{I}_f^0(\nu)$  with inverse $\widetilde{U}\circ S^h$.  This achieves the proof.
\nqed

\subsection{Density of variation processes}

\begin{proposition}
\label{Lemmapn}
Let $(\Omega,\mathcal{A}, (\mathcal{A}_t)_{t\in[0,1]}, \mathcal{P})$ be a complete stochastic basis  and let $p_n$ be the linear operator defined by 
\begin{equation} \label{211} p_n : u\in L^2_a(\mathcal{P},H) \to p_n(u) := \int_0^. \sum_{k=2}^{n-1} n1_{[\frac{k}{n}, \frac{k+1}{n})}(s) \left(u_{\frac{k-1}{n}}- u_{\frac{k-2}{n}} \right) ds \in L^2_a(\mathcal{P},H) \end{equation}
Then for any $u\in L^2_a(\mathcal{P},H)$  $(p_n(u))$ converges to $u$ strongly in $L^2_a(\mathcal{P},H)$. Moreover it satisfies $\mathcal{P}-a.s.$ \begin{equation} \label{linftypn} |p_n(u)|_W \leq 2(n-2)|u|_W \end{equation} and  \begin{equation} \label{linftypn2} |p_n(u)|_H \leq |u|_H \end{equation}
\end{proposition}
\nproof
Inequality~(\ref{linftypn}) follows from the definition. On the other hand~(\ref{linftypn2}) directly follows from Jensen's inequality so that $p_n(u)$ is also a contraction of $L^2(\mathcal{P},H)$ i.e. $|p_n(u)|_{L^2(\mathcal{P},H)} \leq |u|_{L^2(\mathcal{P},H)}$.  For $\epsilon>0$ let  $\beta$ be the primitive of an elementary predicable process such that $$|u- \beta |_{L^2({\mathcal{P}},H)}<   \frac{\epsilon}{2} $$  By primitive of an elementary predicable process  we mean that $\beta$ is of the form \begin{equation} \label{betaform} \beta:=\sum_{k=0}^{n_0} \int_0^. 1_{(t_k,t_{k+1}]}(s) \alpha_k ds \end{equation} where $n_0\in \mathbb{N}$, $(t_k)_{k\in [|0,n_0+1 |]} \subset [0,1]$ is increasing, and where for any $k$, $\alpha_k$ is an element of $L^2_a({\mathcal{P}}, \mathbb{R}^d)$ which is $\mathcal{A}_{t_k}$ measurable.   Such a $\beta$ always exists by a well known result (See \cite{I-W}). Then one can  see that there exists  $N_\beta\in \mathbb{N}$ and a constant $C_\beta \in [0,\infty)$ such that, for any $n>N_\beta$, $$E_{\mathcal{P}}\left[|p_n(\beta)- \beta|_{\mathbb{R}^d}^2 \right] \leq \frac{1}{n} C_\beta(n_0+1) \max_{k\in [|0, n_0+1|]}E_\mathcal{P}\left[|\alpha_k|_{\mathbb{R}^d}^2\right] $$  This shows that $(p_n(\beta))$ converges to $\beta$ for any simple process $\beta$.  Together with the fact that for any $n\in \mathbb{N}$ $p_n$ is a linear contraction this yields 
\begin{eqnarray*}
|p_n(u) -u |_{L^2({\mathcal{P}},H)}   & \leq & |p_n(u-\beta)|_{L^2({\mathcal{P}},H)}  +|p_n(\beta)- \beta |_{L^2({\mathcal{P}},H)}  + | \beta  -u |_{L^2({\mathcal{P}},H)}  \\ & \leq  & 2|u-\beta|_{L^2({\mathcal{P}},H)}  +|p_n(\beta)- \beta |_{L^2({\mathcal{P}},H)}   \\ & \leq & \epsilon + |p_n(\beta)- \beta |_{L^2({\mathcal{P}},H)} 
\end{eqnarray*}
By using the convergence of  $(p_n(\beta))$  to $\beta$ which is a simple process we finally get $$\limsup |p_n(u) -u |_{L^2({\mathcal{P}},H)} \leq \epsilon$$ Since  $\epsilon >0$ is arbitrary, $(p_n(u))$ converges to $u$ strongly in $L^2({\mathcal{P}},H)$.
\nqed

\begin{proposition} 
\label{lemmatangent}
For any $\nu \in \mathcal{P}_W$,  $V_\nu^\infty \cap L^\infty(\nu,H)$ (and then $V_\nu^\infty$ and $V_\nu$) is dense (strongly) in $L^2_a(\nu,H)$.
\end{proposition}
\nproof
First, note that the space $L^\infty(\nu,H)$ is dense in $L^2(\nu,H)$ for the topology of $L^2(\nu,H)$. Indeed or  $u\in L^2(\nu,H)$ by taking \begin{equation} \label{taun} \tau_n:= \inf(\{t : |\pi_tu|_H > n\})\wedge1\end{equation} and $u^n:=\pi_{\tau_n}u:=\int_0^. 1_{[0,\tau_n]}(s) \dot{u}_s ds$, the dominated convergence theorem  yields the convergence of $u^n$ to $u$ in $L^2(\nu,H)$. Hence it is sufficient to prove that $V_\nu \cap L^\infty(\nu,H)$ is dense in $L^\infty(\nu,H)$ for the $L^2(\nu,H)$ topology. To prove this we set $$V^{del}_\nu := \{ p_n(u), u\in L^\infty_a(\nu,H), n\in \mathbb{N} \}$$ where $p_n$ is the operator defined in Proposition~\ref{Lemmapn}.  To prove the density of $V_\nu \cap L^\infty(\nu,H)$ it is sufficient to prove that  \begin{equation} \label{inclm}V^{del}_\nu   \subset V_\nu \cap L^\infty(\nu,H) \end{equation} and that $V^{del}_\nu$ is dense in $L^\infty(\nu,H)$.  First note that by~(\ref{linftypn2}) we already know that $V^{del}_\nu \subset L^\infty(\nu,H)$. Hence the density follows directly by the definition of $V^{del}_\nu$ together with Proposition~\ref{Lemmapn}. So we just have to prove~(\ref{inclm}): we consider a $h\in V^{del}_\nu $ and set \begin{equation} \label{uhdefp} U^h:= U +  h\end{equation} where $U\in \mathcal{I}_f^0(\nu)$. We consider two measurable mappings $X:\Omega \to W$ and $Y :\Omega \to W$ defined on the same probability space $(\Omega,\mathcal{A},\mathcal{P})$ such that $X_\star \mathcal{P} = Y_\star \mathcal{P}=\nu$ and we assume that \begin{equation} \label{thetahequal} U^h(X)= U^h(Y)\end{equation} By Proposition~\ref{propisom}, to obtain that $U^h \in \mathcal{I}_f^0(\nu)$ it is sufficient to prove that we necessarily have $\mathcal{P}-a.s.$ $X=Y$, and that for all $t\in [0,1)$, \begin{equation} \label{mpfmjc}\mathcal{G}_t^{U^h} =\mathcal{F}_t^\nu\end{equation} We postpone~(\ref{mpfmjc}) to the end of the proof and we set $$\tau:= \inf(\{ t : X_t \neq Y_t\})\wedge 1$$ Since $U\in \mathcal{I}_f^0(\nu)$ we also have \begin{equation} \label{isompf} \tau = \inf\left(\left\{ t : U_t(X) \neq U_t(Y)\right\}\right)\wedge1\end{equation} Note that,  since $U\in \mathcal{I}_f^0(\nu)$ there exists a $\mathcal{B}(\mathbb{R}^d)^{{U_0}_\star \eta}/ \mathcal{B}(\mathbb{R}^d)-$measurable mapping $f : \mathbb{R}^d \to \mathbb{R}^d$ so that $\nu-a.s.$ $$W_0= f(U_0) $$ by~(\ref{thetahequal}), since $\nu-a.s.$ $h_0=h_1=0$ it yields $\mathcal{P}-a.s.$ \begin{equation} \label{val0in} X_0= Y_0 \end{equation} On the other hand (see Proposition~\ref{Lemmapn}) there exists a $\lambda>0$ such that $h$ is adapted to the filtration $(\mathcal{H}_t^\lambda)$ where for $t\geq \lambda$ (resp. $t<\lambda$) $\mathcal{H}_t^\lambda= \mathcal{F}^\nu_{t-\lambda}$ (resp. $\mathcal{H}_t= \mathcal{B}^0(W)^\nu$). Using this together with~(\ref{val0in}) we obtain $\mathcal{P}-a.s.$ $$\inf\{ t : h_t(X) \neq h_t(Y)\} \geq (\tau + \lambda)\wedge 1$$ together with~(\ref{thetahequal}) and~(\ref{isompf}) this latter inequality reads $$\tau \geq (\tau + 
\lambda) \wedge 1$$ so that $\mathcal{P}-a.s.$ $\tau=1$ and $X= Y$. Hence to show that $U^h\in \mathcal{I}_f^0$ it is now sufficient to prove~(\ref{mpfmjc}) for all $t\in[0,1)$. The first inclusion is trivial, and we just have to prove that for any $t\in[0,1)$ \begin{equation} \label{includfilt} \mathcal{F}_t^\nu \subset \mathcal{G}^{U^h}_t\end{equation} We choose $N_\lambda\in \mathbb{N}$ such that $\lambda N_\lambda >2$ and we set $t_i= \frac{i}{N_\lambda}$ for a $i\in[0, N_\lambda-1]\cap\mathbb{N}$. For $t\in [t_0,t_1)$  we have  $\nu-a.s.$ $h_t=0$ for all $t\in [t_0,t_1)$ so that $ \sigma(U_s, s\leq t)^\nu = \sigma(U^h_s, s\leq t)^\nu$. Since $U\in \mathcal{I}_f^0(\nu)$  we have $\mathcal{G}^{U^h}_t= \mathcal{G}^{U}_t= \mathcal{F}_t^\nu$ for $t\in[t_0,t_1)$. 
Finally we assume that~(\ref{includfilt}) holds for any $t\in[t_0,t_i)$ and we take $t\in [t_i,t_{i+1})$. Almost surely with respect to $\nu$, we have $U_t = U^h_t - h_t$. But $h_t$ is  $\mathcal{F}^\nu_{t-\lambda}$ measurable and therefore $
\mathcal{F}^\nu_{t_{i-1}}$ measurable so that by the induction hypothesis it is also $\mathcal{G}^{U^h}_{t_{i-1}}$ 
measurable and hence $\mathcal{G}^{U^h}_t-$measurable. Thus, for all $t\in[t_i, t_{i+1})$  \begin{equation} \label{mfnxj} \sigma(U_s, s
\leq t) \subset \mathcal{G}_t^{U^h}\end{equation} Using the fact that $U\in \mathcal{I}_f^0(\nu)$ together with the right continuity of $(\mathcal{G}_t^{U^h})$,~(\ref{mfnxj}) yields   $$\mathcal{F}_t^\nu = 
\mathcal{G}^{U}_t\subset  \mathcal{G}_t^{U^h}$$ for all $t\in [t_i,t_{i+1})$, and by induction for all $t\in[0,1)$. This proves that $h\in V_\nu \cap L^\infty
(\nu,H)$ so that ~(\ref{inclm}) holds. The proof is complete.
\nqed

\subsection{Variation processes with vanishing endpoints}

\begin{proposition}
\label{densepropl}
For any $\nu \in \mathcal{P}_W$ the space $L^\infty(\nu,W) \cap L^2_a(\nu, H_{0,0})$ is dense in $L^2_a(\nu, H_{0,0})$ for the topology of $L^2(\nu,H).$ \end{proposition}
\nproof
For $u:=\int_0^. \dot{u}_s ds \in L^2_a(\nu, H_{0,0})$ and $\tau$ a $(\mathcal{F}_t^\nu)-$stopping time we set
\begin{equation} \label{ktaudef} k^\tau[u]:= \int_0^{.\wedge \tau} \dot{u}_s ds   -\int_0^. \frac{1_{(\tau, 1)}(s)}{1-\tau} u_\tau ds \end{equation}
By Jensen's inequality~(\ref{ktaudef}) defines a linear and continuous operator $k^\tau : L^2_a(\nu, H_{0,0}) \to  L^2_a(\nu, H_{0,0})$. More precisely for $u\in L^2_a(\nu, H_{0,0})$, since $\nu-a.s.$ $u_1=0$, Jensen's inequality yields  $\nu-a.s.$ \begin{equation}\label{jktau1} |k^{\tau}[u]|_H \leq |u|_H \end{equation} and   \begin{equation}\label{jktau2} |k^{\tau}[u]|_W \leq 2 \sup_{t\leq\tau } |u_t| = 2|\pi_\tau u|_W \end{equation}  Since $|.|_W\leq |.|_H$ by~(\ref{jktau2}) we have, $\nu-a.s.,$  \begin{equation} \label{jktau3}  |k^{\tau}[u]|_W \le 2 |\pi_\tau u |_H \end{equation}  For any $n\in \mathbb{N}$, let $(\tau_n)$ be the sequence of stopping time associated to $u$ by~(\ref{taun}).  We define a sequence $(u^n)$ by setting $u^n:= k^{\tau_n}[u]$. For all $n\in \mathbb{N}$ $u^n\in L^2_a(\nu, H_{0,0})$ and by~(\ref{jktau3}) $u^n \in L^\infty(\nu, W)$. Hence $(u^n)_{n\in \mathbb{N}} \subset L^2_a(\nu, H_{0,0}) \cap L^\infty(\nu,W)$. Finally by~(\ref{jktau1}) the dominated convergence theorem yields the convergence of $(u^n)$ to $u$ and therefore the density of $L^\infty(\nu,W) \cap L^2_a(\nu, H_{0,0})$. 
\nqed

\begin{lemma}
\label{denselemma}
For any $\nu \in \mathcal{P}_W$, $V_\nu^{0,\infty}$ is dense in $L^2_a(\nu, H_{0,0})$.
\end{lemma}
\nproof
For $n\in \mathbb{N}$, $n\geq 2$ we set  \begin{equation} \label{qndef}q_n :  u \in L^2_a(\nu,H) \to q_n(u):=\int_0^. 1_{[1-\frac{1}{n},1]}(s) n u_{1-\frac{2}{n}} ds \in L^2_a(\nu,H) \end{equation} and  \begin{equation} \label{rndef} r_n :=  p_n - q_n \end{equation} where $p_n$ is defined in Proposition~\ref{Lemmapn}. In particular for any $u\in L^2_a(\nu,H)$ $\nu-a.s.$ \begin{equation} \label{eqnW}|q_n(u)|_W \leq |u|_W \end{equation} We set $$V^{0,del}_\nu := \{ r_n(u) : u\in L^2_a(\nu,H_{0,0}) \cap L^\infty(\nu,W) , n\in \mathbb{N} \}$$ By Proposition~\ref{densepropl} it is sufficient to prove that $V^{0,del}_\nu$ is a dense subset of $L^\infty(\nu,W) \cap L^2_a(\nu, H_{0,0})$ for the $L^2(\nu,H)$ topology and that \begin{equation} \label{tdelincl}  V^{0,del}_\nu  \subset V_\nu^{0,\infty} \end{equation} The density follows by Jensen's inequality. Indeed for $u\in L^\infty(\nu,W) \cap L^2_a(\nu, H_{0,0})$  since $\nu-a.s.$ $u_1=0$, Jensen's inequality yields $\nu-a.s.$, $$|q_n(u)|_H^2 \leq 2\int_{1-\frac{2}{n}}^1 |\dot{u}_s|^2 ds $$ and  $$|q_n[u]|_{L^2(\nu,H)}\leq \sqrt{2} | (I_H- \pi_{1-\frac{2}{n}}) u|_{L^2(\nu,H)}$$ Bby the dominated convergence theorem $q_n(u)$ converges strongly to $0_{L^2(\nu,H)}$ in $L^2(\nu,H)$. Hence by Proposition~\ref{Lemmapn}, $(r_n(u))$ converges to $u$. On the other hand by~(\ref{linftypn}) and~(\ref{eqnW}) \begin{equation} \label{odi1}V^{0,del}_\nu \subset L^\infty(\nu,W) \end{equation} and by~(\ref{211}),~(\ref{qndef}) and~(\ref{rndef}) we have $\nu-a.s.$ $r_n(1)=0$ so that we also have \begin{equation} \label{odi2} V^{0,del}_\nu \subset L^2_a(\nu, H_{0,0}) \end{equation} Finally for a $U\in \mathcal{I}_f^0(\nu)$ and $h\in V^{0,del}_\nu$, similarly to the proof of  Proposition~\ref{lemmatangent}  we obtain that $U+h \in \mathcal{I}_f^0(\nu)$ i.e. $V^{0,del}_\nu\subset V^\nu$. Together with~(\ref{odi1}) and~(\ref{odi2}) we derive ~(\ref{tdelincl}), by which the result follows.
\nqed

\section{Transformations of $\mathbb{S}$}
\label{section3}

\subsection{Transformation formulas}
The aim of this subsection is to prove Proposition~\ref{lemmaforback2} and Lemma~\ref{caspartp9} which are the only results of this subsection that will be used in the sequel. Before we prove Proposition~\ref{lemmaforback}  in order to derive Proposition~\ref{propframe}. Both results  can be skipped; however this latter proposition justifies the definitions of the set $V_\nu$. Despite their apparent generality Proposition~\ref{lemmaforback}  and Proposition~\ref{propframe} assume some further integrability assumptions due to the fact that their proofs involve the dual predicable projection. In the particular case of Lemma~\ref{caspartp9} and Proposition~\ref{lemmaforback2} we stress that the use of isomorphisms of filtered probability spaces allows us to drop this condition.

\begin{proposition}
\label{lemmaforback}
Given a complete stochastic basis $(\Omega,\mathcal{A}, (\mathcal{A}_t)_{t\in[0,1]}, \mathcal{P})$, let $U:\Omega \to W$ be a $\mathcal{A}/\mathcal{B}(W)$-measurable mapping such that $\mathcal{P}$-a.s. for all $t\in[0,1]$ \begin{equation} \label{udelmf} U_t = U_0 + M^u_t + \int_0^t \dot{u}_s ds  \end{equation} where $(t,\omega) \in [0,1]\times \Omega \to \dot{u}_t\in \mathbb{R}^d$ is $(\mathcal{A}_t)-$predicable with, for all $T<1$, \begin{equation} \label{dominconv}  E_{\mathcal{P}}\left[\int_0^T |\dot{u}_s|_{\mathbb{R}^d} ds \right] <\infty \end{equation} and where $(M^u_t)$ is a continuous $\mathbb{R}^d-$valued $(\mathcal{A}_t)-$local martingale such that $\mathcal{P}-a.s.$ for all $i,j\in [1,d]$ \begin{equation} \label{mucoch} <(M^u)^i, (M^u)^j>=\int_0^. (\alpha_s^u)^{ij} ds  \end{equation} for some $\mathcal{M}_d(\mathbb{R})-$valued $(\mathcal{A}_t)-$predicable process $(\alpha_s^u)$. By setting  $$\nu:=U_\star\mathcal{P}$$ we have $\nu \in \mathbb{S}$. Moreover $\mathcal{P}-a.s.$  \begin{equation} \label{crochetmp}\int_0^.  \alpha_t^\nu \circ U dt =  \int_0^. \alpha^u_t dt \end{equation}and $\mathcal{P}-a.s.$ \begin{equation}\label{veqlfb} \int_0^. v^{\nu}_t\circ U dt =  \int_0^. E_{\mathcal{P}}[\dot{u}_t| \mathcal{G}^U_t] dt  \end{equation} where $(E_{\mathcal{P}}\left[\dot{u}_t \middle| \mathcal{G}_t^U\right])_{t\in[0,1]}$ denotes the optional projection of $(\dot{u}_t)$ on $(\mathcal{G}_t^U)_{t\in[0,1]}$, the usual augmentation of the filtration generated by $(U_t)_{t\in[0,1]}$.
\end{proposition}
\nproof 
Let $(\widehat{u}_t)$ be the dual predicable projection of $(u_t)$ (whose variations are locally integrable by~(\ref{dominconv})) on $(\mathcal{G}_t^U)$. In particular \begin{equation} \label{mdpme} \widehat{u} = \int_0^.E_{\mathcal{P}}\left[\dot{u}_t \middle| \mathcal{G}_t^U\right] dt \end{equation}  Then by setting $\mathcal{P}-a.s.$ for all $t\in[0,1]$ \begin{equation} \label{sddlfq} \widehat{M}^u_t:= M^u_t + u_t -\widehat{u}_t  \end{equation} we have that $\mathcal{P}-a.s.$ for all $t\in[0,1]$ \begin{equation} \label{Mwhu} U_t - U_0 - \int_0^t E_{\mathcal{P}}\left[\dot{u}_s \middle| \mathcal{G}_s^U\right] ds = \widehat{M}^u_t\end{equation} so that  $(\widehat{M}^u_t)$ is $(\mathcal{G}_t^U)-$adapted. Let $b$  (resp. $M$) be the $U_\star\mathcal{P}-$equivalence class of $\mathcal{B}(W)^{U_\star\mathcal{P}}/\mathcal{B}(W)$-measurable mappings $(\mathcal{F}_t^{U_\star\mathcal{P}})-$ adapted such that $\mathcal{P}-a.s.$ \begin{equation} \label{v} b \circ U =  \widehat{u}\end{equation} resp. such that $\mathcal{P}-a.s.$ \begin{equation} \label{M} M\circ U= \widehat{M}^u \end{equation} whose existence is insured by Proposition~\ref{propadapt}.  Since $\mathcal{P}-a.s.$ $\widehat{u}$ is absolutely continuous, by~(\ref{v}), $U_\star\mathcal{P}-a.s.$ $b$ is absolutely continuous.  Therefore, there exists  a $(\mathcal{F}_t^{U_\star\mathcal{P}})-$predicable process $(v_t)$ on $(W,\mathcal{B}(W)^{U_\star\mathcal{P}}, U_\star\mathcal{P})$ so that $\mathcal{P}-a.s. $ $$b = \int_0^. v_s ds $$ and by~(\ref{dominconv}) and~(\ref{v}) we have, for all $T<1$, $$E_{U_\star\eta}\left[\int_0^T |v_t|_{\mathbb{R}^d} dt \right] < \infty$$  Note that the $(\mathcal{G}_t^U)-$local martingale $(\widehat{M}^u_t)$ is reduced by the sequence of $(\mathcal{G}_t^U)$ stopping times \begin{equation} \label{whtnd}  \widehat{\tau}_n := \inf(\{ t\in [0,1] :  |\widehat{M}^u_t| > n\}) \end{equation}   Indeed by the stopping theorem and the dominated convergence theorem~(\ref{dominconv}) easily implies that the process $(M^u_t)$ is reduced to $(\mathcal{A}_t)-$martingale by the sequence $(\widehat{\tau}_n)_{n\in \mathbb{N}}$, while by~(\ref{M}) and~(\ref{whtnd}) $(\widehat{\tau}_n)$ is also a family of $(\mathcal{G}_t^U)$ stopping times.  Therefore together with~(\ref{sddlfq}) and using the inclusion $(\mathcal{G}_.^U)\subset \mathcal{A}_.$  we obtain that $(\widehat{\tau}_n)$ also reduce $(\widehat{M}^u_t)$ to  $(\mathcal{G}_t^U)$ martingales.  By Dynkin's Lemma  and since $\omega \to \widehat{\tau}_n(\omega)$ is $\mathcal{G}_1^U$ measurable, for $n\in  \mathbb{N}$ we denote by $\tau_n$ the $(\mathcal{F}_t^{U_\star \mathcal{P}})-$stopping time such that $\mathcal{P}-a.s.$ $$\widehat{\tau}_n = \tau_n\circ U$$  Since $\widehat{M}^u$ is a $(\mathcal{G}_t^U)-$local martingale reduced by $\widehat{\tau}_n$, it is straightforward to check that $M$ is a $(\mathcal{F}_t^{U_\star\mathcal{P}})$ local martingale under $U_\star\mathcal{P}$ (reduced by $(\tau_n)$). Indeed for $n\in \mathbb{N}$ and  $s\leq t$ we have $\mathcal{P}-$a.s. $$E_{\mathcal{P}}\left[\widehat{M}^u_{t\wedge \widehat{\tau}_n} | \mathcal{G}_s^U\right] = E_{\mathcal{P}}\left[(M_{t\wedge \tau_n}) \circ U \middle| \mathcal{G}_s^U\right] = E_{U_\star\mathcal{P}}\left[M_{t\wedge \tau_n} \middle| \mathcal{F}_s^{U_\star\mathcal{P}}\right]\circ U$$ so that $\mathcal{P}-a.s.$  $$ E_{U_\star\mathcal{P}}\left[M_{t\wedge \tau_n} \middle| \mathcal{F}_s^{U_\star\mathcal{P}}\right]\circ U = E_{\mathcal{P}}\left[\widehat{M}^u_{t\wedge \widehat{\tau}_n} \middle| \mathcal{G}_s^U\right] = \widehat{M}^u_{s\wedge\widehat{\tau}_n} = M_{s\wedge \tau_n}\circ U $$  and $U_\star\mathcal{P}-a.s.$ $$E_{U_\star\mathcal{P}}\left[M_{t\wedge \tau_n} \middle| \mathcal{F}_s^{U_\star\mathcal{P}}\right] = M_{s\wedge \tau_n}$$ By writing the Dolean's approximations of the predicable quadratic co-variation process as a limit of finite sums we obtain  $\mathcal{P}-a.s.$ for all $t\in [0,1]$, \begin{equation} \label{eqcrol}<M^i, M^j>_t\circ ~U = <(M^u)^i,(M^u)^j >_t\end{equation} since the process of the right hand term is absolutely continuous, the process of the left hand term too and we have the existence of a $(\mathcal{F}_t^{U_\star\mathcal{P}})-$predicable process $(\alpha_t)$  on the probability space $(W,\mathcal{F}_t^{U_\star\mathcal{P}}, U_\star\mathcal{P})$ such that $U_\star\mathcal{P}-a.s.,$ $$<M^i, M^j>= \int_0^. (\alpha_s)^{ij} ds$$ Finally by reporting~(\ref{v}) and~(\ref{M}) in~(\ref{Mwhu}) we obtain, $U_\star\mathcal{P}-a.s.$ for all $t\in[0,1]$, $$W_t - W_0 - \int_0^t v_s ds = M_t$$ by which the result follows.
\nqed

\begin{proposition}
\label{propframe}
For $\nu \in \mathbb{S}$  whose characteristics are denoted $(\int_0^. v_t^\nu dt, \int_0^. \alpha^\nu_t dt)$ (see Definition~\ref{defvelocity}), let $\nu-a.s.$  $h:=\int_0^. \dot{h}_s ds \in L^0_a(\nu,W_{abs})$ where $(\dot{h}_s)$ is a $(\mathcal{F}_t^\nu)$-predicable process and let $(\theta_s)$ be a $\mathcal{M}_d(\mathbb{R})$-valued $(\mathcal{F}_t^\nu)$-predicable process on  $(W,\mathcal{B}(W)^\nu,\nu)$ such that for all $T<1$ \begin{equation} E_\nu\left[\int_0^T  |\theta_s v^\nu_s +\dot{h}_s |_{\mathbb{R}^d} ds\right] <\infty \end{equation} and $\nu-a.s.$ for all $i,j$  \begin{equation} \label{csmarp} \int_0^1 (\theta_s^{i,j})^2 {\alpha_s^\nu}^{jj} ds <\infty \end{equation} Furthermore define $R\in L^0(\nu,W)$ by $\nu-a.s.$ for all $t\in[0,1]$ \begin{equation} \label{defdeR}  R_t = f(W_0) + \int_0^t \theta_s dW_s \end{equation}  where $f :\mathbb{R}^d\to \mathbb{R}^d$ is any $\mathcal{B}(\mathbb{R}^d)^{W_0\star\nu}/ \mathcal{B}(W)-$measurable function, and $U\in L^0(\nu,W)$ by  \begin{equation} \label{defdeU} U = R+ h \end{equation}   i.e. $\nu-a.s.$ for all $t\in[0,1]$  \begin{equation} \label{359'} U_t:= f(W_0) + \int_0^t \theta_s dW_s + \int_0^t \dot{h}_s ds  \end{equation} Then $U_\star \nu \in \mathbb{S}$ and we have, $\nu-a.s.$,
\begin{equation} \label{vfpour}  \int_0^. v_s^{U_\star \nu}\circ U  ds =  \int_0^. E_\nu\left[\theta_s v_s^{ \nu} +\dot{h}_s \middle|\mathcal{G}_s^U\right]   ds      \end{equation} and $\nu-a.s.$   \begin{equation}  \label{afpour} \int_0^. \alpha^{U_\star \nu}_s\circ U ds  =  \int_0^. \theta_s \alpha_s^\nu \theta_s^{\dagger} ds \end{equation} where $\left(E_\nu\left[\theta_t v_t^{ \nu} +\dot{h}_t \middle|\mathcal{G}_t^U\right]\right)_{t\in[0,1]}$ denotes the optional projection of the process $$(\theta_t v_t^{ \nu} +\dot{h}_t)_{t\in[0,1]}$$ on the filtration $(\mathcal{G}_t^U)$. In particular if $U\in \mathcal{I}_f^0(\nu)$, for all $\epsilon \in \mathbb{R}$
$\nu-a.s.$ \begin{equation} \label{365M} \int_0^. v_s^{U_\star \nu}\circ U   ds  =  \int_0^. \left(\theta_s v_s^{\nu } +\epsilon \dot{h}_s\right)   ds     \end{equation}
\end{proposition}
\nproof
Denote $M^\nu$ the martingale part of $\nu$. Then $\nu-a.s.$ for all $t\in[0,1]$, \begin{equation} \label{373m} R_t= f(W_0) +  M^u_t + \int_0^t \theta_s v_s^\nu ds \end{equation}  and  $$ U_t= f(W_0) + M^u_t  + \int_0^t ( \theta_s v_s^\nu +\dot{h}_s ) ds $$where $(M_t^u)$ is the process defined on $(W, \mathcal{B}(W)^\nu, \nu)$ by $\nu-a.s.$ for all $t\in[0,1]$ \begin{equation} \label{mun25} M_t^u := \int_0^t \theta_s dM_s^\nu \end{equation} By~(\ref{csmarp}) $(M_t^u)_{t\in[0,1]}$ is a continuous $(\mathcal{F}_t^\nu)$-local martingale on $(W,\mathcal{B}(W)^\nu, \nu)$  and $\nu-a.s.$ for all $t\in[0,1]$ \begin{equation} \label{mun26} < (\int_0^. \theta_s dM_s^\nu)^i, (\int_0^. \theta_s dM_s^\nu)^j>_t=\sum_{l,m}\int_0^t \theta_s^{il} \theta_s^{jm} (\alpha^\nu_s)^{lm}dt\end{equation} for $i,j\in[1,d]$.  Hence by applying Proposition~\ref{lemmaforback} to $U$ on $(W,\mathcal{B}(W)^\nu,\nu)$ with the filtration $(\mathcal{F}_t^\nu)$ we obtain $U_\star \nu \in \mathbb{S}$ and~(\ref{vfpour}) and~(\ref{afpour}). The end of the claim follows from the definition of $\mathcal{I}_f^0(\nu)$.
\nqed

\begin{proposition}
\label{lemmaforback2}
For $\nu \in \mathbb{S}$, let $\nu-a.s.$ $k^u:=\int_0^. \dot{u}_s ds \in L^0_a(\nu,W_{abs})$ and let $(M^u_t)$ be a continuous $(\mathcal{F}_t^\nu)$ local martingale on $(W,\mathcal{B}(W)^\nu,\nu)$ vanishing at $t=0$, such that $$<(M^u)^i, (M^u)^j> = \int_0^. (\alpha_s^u)^{i,j} ds $$ for some $\mathbb{R}^d$ (resp. $\mathcal{M}_d(\mathbb{R})$) valued $(\mathcal{F}_t^\nu)$ predicable process $(\dot{u}_s)$ (resp. $(\alpha^u_s)$) on $(W, \mathcal{B}(W)^\nu,\nu)$.  Moreover consider a continuous process defined by $\nu-a.s.$ for all $t\in[0,1]$  \begin{equation} \label{UdefcasM} U_t := f(W_0) + M^u_t + k^u_t\end{equation} where $f : \mathbb{R}^d\to \mathbb{R}^d$ is as in Proposition~\ref{propframe}. Assume  that $U\in \mathcal{I}_f^0(\nu)$. Then $$ U_\star \nu \in \mathbb{S}$$ Moreover $\nu-a.s.$  \begin{equation}\label{1R311M} \int_0^. v_s^{U_\star\nu}\circ U ds   = \int_0^. \dot{u}_s ds \end{equation}  and $\nu-a.s.$ \begin{equation} \label{2R311M} \int_0^. \alpha_s^{U_\star\nu}\circ U ds   = \int_0^. \alpha_s^u ds \end{equation}  
\end{proposition}
\nproof
 We denote $\widetilde{U}\in \mathcal{I}_f^0(U_\star\nu)$ the inverse of $U$ and we set $\eta:=U_\star \nu$. Using this we denote $\widetilde{M}$ (resp. $\widetilde{b}$) the elements of $L^0_a(\eta, W)$ (resp. of $L^0_a(\eta, W_{abs})$) defined by $\eta-a.s.$ \begin{equation} \label{btildedefemM} \widetilde{M}:= M^u\circ \widetilde{U}\end{equation} and by $\eta-a.s.$ \begin{equation} \label{btildefemM2} \widetilde{b}:=  k^u \circ \widetilde{U} \end{equation} Denote $(\tau_n)$ (resp. $(\widetilde{\tau}_n)$) a sequence of $(\mathcal{F}_t^\nu)$ stopping times reducing $M^u$ to a $(\mathcal{F}_t^\nu)$-martingale (resp., since $\mathcal{G}_.^{\widetilde{U}} =\mathcal{F}_.^\eta$, of $(\mathcal{F}_t^\eta)$-stopping times defined for all $n\in \mathbb{N}$ by $\eta-a.s. $ $\widetilde{\tau}_n:= \tau_n\circ \widetilde{U}$), similarly to the proof of Proposition~\ref{lemmaforback} and using both that $U$ preserves the filtration and that $\widetilde{U}_\star U_\star \nu= \nu$ it is straightforward to obtain that, for all $n\in \mathbb{N}$, $(\widetilde{\tau}_n)$ reduces $\widetilde{M}$ to a $(\mathcal{F}_t^\eta)$ martingale on $(W,\mathcal{B}(W)^\eta, \eta)$. Since $U\in \mathcal{I}_f^0(\nu)$ by~(\ref{UdefcasM}) $\eta-a.s.$, \begin{equation} \label{W0UT0} W_0 = U_0\circ \widetilde{U}= f(\widetilde{U}_0) \end{equation}  By~(\ref{UdefcasM}),~(\ref{btildedefemM}),~(\ref{btildefemM2}) and~(\ref{W0UT0}) since $\eta-a.s.$ $U\circ\widetilde{U}= I_W$ we obtain $\eta-a.s.$ for all $t\in[0,1]$, \begin{equation} \label{coordlmpM} W_t = W_0 + \widetilde{M}_t  + \widetilde{b}_t \end{equation} so that $\eta\in \mathbb{S}$. Then the result follows similarly to the proof of Proposition~\ref{lemmaforback} by the hypothesis $U\in \mathcal{I}_f^0(\nu)$. Indeed by~(\ref{coordlmpM}) and~(\ref{btildefemM2}) we obtain $\nu-a.s.$ $b^\eta\circ U= k^u$ which yields~(\ref{1R311M}) while~(\ref{2R311M}) follows similarly.
\nqed

\begin{lemma}
\label{caspartp9}
For $\nu \in \mathbb{S}$  whose characteristics are denoted $(\int_0^. v_t^\nu dt, \int_0^. \alpha^\nu_t dt)$ (see Definition~\ref{defvelocity}), let $\nu-a.s.$ $h:=\int_0^. \dot{h}_s ds \in L^0_a(\nu,W_{abs})$ and let $(\theta_s)$ be a $\mathcal{M}_d(\mathbb{R})$-valued predicable process on $(W,\mathcal{B}(W)^\nu,\nu)$ such that $\nu-a.s.$~(\ref{csmarp}) holds for all $i,j\in [1,d]$ and  $\nu-a.s.$, \begin{equation} \int_0^1  |\theta_s v^\nu_s |_{\mathbb{R}^d} ds <\infty \end{equation}  For $f : \mathbb{R}^d \to \mathbb{R}^d$ as in Proposition~\ref{propframe}, let $R\in L^0(\nu,W)$ be given by~(\ref{defdeR}) and  $U\in L^0(\nu,W)$ by~(\ref{defdeU}). If we further assume that $U\in \mathcal{I}_f^0(\nu)$ then  $U_\star \nu \in \mathbb{S}$ and $\nu-a.s.$ we have~(\ref{afpour}) and~(\ref{365M}).  In particular for $h\in V_\nu$, denoting  $\nu-a.s.$ $$\tau_h := I_W + h$$ where $I_W$ denotes the $\nu-$equivalence class of mappings $\nu-a.s.$ equal to the identity map, we have for all $\epsilon\in \mathbb{R}$, ${\tau_{\epsilon h}}_\star \nu \in \mathbb{S}$. Moreover $\nu-a.s.$, \begin{equation} \label{1R311} \int_0^. v_s^{{\tau_{\epsilon h}}_\star \nu}\circ \tau_{\epsilon h}   ds =  \int_0^. (v_s^{ \nu} + \epsilon \dot{h}_s ) ds   \end{equation} and \begin{equation} \label{2R311} \int_0^. \alpha^{{\tau_{\epsilon h}}_\star \nu}_s\circ \tau_{\epsilon h}  ds =\int_0^. \alpha^{\nu}_s ds \end{equation}
\end{lemma}
\nproof
Denote $M^\nu$ the martingale part of $\nu$. Then $\nu-a.s.$ for all $t\in[0,1]$ $R_t$ is given by~(\ref{373m}) and $U_t$ by \begin{equation} \label{Udefcas} U_t= f(W_0) + M^u_t  + k^u_t \end{equation} where $$k^u := \int_0^. ( \theta_s v_s^\nu +\dot{h}_s ) ds$$ and where $(M_t^u)$ is the process defined on $(W, \mathcal{B}(W)^\nu, \nu)$ by~(\ref{mun25}). By~(\ref{csmarp}), $(M_t^u)_{t\in[0,1]}$ is a continuous $(\mathcal{F}_t^\nu)$-local martingale on $(W,\mathcal{B}(W)^\nu, \nu)$  and $\nu-a.s.$ we also have~(\ref{mun26}). Since $U\in \mathcal{I}_f^0(\nu)$  the result follows by Proposition~\ref{lemmaforback2}. Moreover since $I_W\in \mathcal{I}_f^0(\nu)$ by definition of $V_\nu$, $\tau_h\in V_\nu$ so that~(\ref{1R311}) and~(\ref{2R311}) follows as a particular case.
\nqed

\subsection{Lift of transformations on the space to transformations of $\mathbb{S}$}

The next proposition is general, however it is formulated to provide an insight on the behaviour of transformations which are close to the identity
\begin{lemma}
\label{constantelemm}
Let $h : (t,x)\in [0,1] \times \mathbb{R}^d \to h(t, x)\in \mathbb{R}^d$ be a $C^{1,2}$ function ($C^1$ in $t$, $C^2$ in $x$), and set $u(t,x) := h(t,x)-x$. Assume also that, for any $t\in [0,1]$, $$h_t:  x\in \mathbb{R}^d \to h_t(x):= h(t,x) \in \mathbb{R}^d$$ is an homeomorphism, whose inverse $j_t $ is such that $(t,x)\in [0,1]\times \mathbb{R}^d \to j_t(x) \in \mathbb{R}^d$ is continuous.  Denote $G: W\rightarrow W$ (resp. $\widetilde{G}: W\rightarrow W$) the mapping defined for all $(t,\omega)\in [0,1]\times W$  by 
$$G(\omega)(t):= h(t, \omega(t)) =\omega(t) +  u(t,\omega(t))$$  (resp. $\widetilde{G}(\omega)(t):= j_t(\omega(t))$).  For any $\nu \in \mathbb{S}$, $G$  induces an element $
\Gamma$ of $\mathcal{I}_f^0(\nu)$ ($\Gamma$ is the element of $L^0(\nu,W)$ such that $\nu-a.s.$ $\Gamma= G$) whose inverse $\widetilde{\Gamma}\in \mathcal{I}_f^0(\Gamma_\star \nu)$  is induced by $\widetilde
{G}$.  Then, for all $\nu \in \mathbb{S}$ we have $${\Gamma}_\star\nu \in \mathbb{S}$$ Moreover $ \nu-a.s.$  \begin{eqnarray*}  \int_0^. v_t^{{\Gamma}_\star\nu}\circ \Gamma dt& = &  \int_0^. \left(\partial_t h(t,W_t) + (v_t^\nu. \nabla) h(t,W_t) + 
\sum_{i,j}\frac{{\alpha^\nu_t}^{i,j}}{2} \partial^2_{i,j}h(t,W_t) \right) dt \\ & = &    \int_0^. \left(v_t^\nu + \partial_t u_t(W_t) + (v_t^\nu. \nabla) u_t(W_t) + \sum_{i,j}\frac{{\alpha^\nu_t}^
{i,j}}{2} \partial^2_{i,j}u(t,W_t) \right) dt \end{eqnarray*} and  $\nu-a.s.$   \begin{eqnarray*}  \int_0^. \alpha_t^{{\Gamma}_\star\nu}\circ \Gamma  dt &  =  &  \int_0^. \left(( \nabla 
h)(t,W_t).(\alpha_t^{\nu}). ( (\nabla h)^{\dagger})(t,W_t)\right) dt 
\\ & = & \int_0^. \left( \alpha_t^\nu +  (\nabla u_t)(W_t).\alpha_t^{\nu}+  \alpha_t^{\nu}. (\nabla u_t)^{\dagger})(W_t) +  (\nabla u_t)(W_t).\alpha_t^{\nu}.(\nabla u_t)^{\dagger}(W_t)  \right) dt  \end{eqnarray*} where $(\nabla h)^{i,j}(x)= \frac{\partial h^i}{\partial x^j}(t,x)$ (similarly for $u$) \end{lemma}
\nproof
By definition for all $t\in[0,1]$ and for all $x\in \mathbb{R}^d$ $$j_t(h_t(x))=h_t(j_t(x))=x$$  Hence  for all $\omega \in W$ $$\widetilde{G}\circ {G} (\omega)(t)= j_t(G_t(\omega))= j_t(h_t(W_t(\omega))= \omega(t)$$ so that $G$ (resp. $\widetilde{G}$) induces an isomorphism of probability spaces (resp. its inverse). On the other hand both ${\Gamma}$ and $\widetilde{\Gamma}$ are adapted to the respective canonical filtrations. This proves that $\Gamma\in \mathcal{I}_f^0(\nu)$. Set $$k^u:= \int_0^. \left(\partial_t h(t,W_t) + (v_t^\nu. \nabla) h(t,W_t) + \sum_{i,j}\frac{{\alpha^\nu_t}^{i,j}}{2} \partial^2_{i,j}h(t,W_t) \right)dt$$ and $$M^u:= \sum_{i=1}^d \left(\int_0^. \partial_{j}h^i(t,W_t)  d(M_t^\nu)^j \right) e_i$$ where $(e_i)_{i=1,d}$ is the canonical orthogonal basis of $\mathbb{R}^d$ and set $\mathcal{P}-a.s.$ for all $t\in[0,1]$ $$U_t:= h_0(W_0) + M^u_t  + k^u_t$$    On the other hand $(M_t^u)$ is a $({\mathcal{F}_t}^\nu)$ local martingale and $$<(M^u)^i, (M^u)^j>= \int_0^.\sum_{l,m} \partial_{l}h^i(t,W_t)  \partial_{m}h^j(t,W_t) (\alpha_t^{\nu})^{l,m} dt $$ By It\^{o}'s formula, $\nu-a.s.$ for all $t\in[0,1]$  $$W_t\circ \Gamma = U_t$$ Therefore the result follows from Proposition~\ref{lemmaforback2}  \nqed

\begin{definition}
\label{lift}
 For $\nu \in \mathbb{S}$ and  $h : (t,x)\in [0,1] \times \mathbb{R}^d \to h_t(x)\in \mathbb{R}^d$ which satisfy the hypothesis of Lemma~\ref{constantelemm} we call the isomorphism of filtered probability spaces $\Gamma\in \mathcal{I}_f^0(\nu)$ associated to $h$ by Lemma~\ref{constantelemm} the\textbf{ lift }of $h$ on $\mathbb{S}$ at $\nu$. 
\end{definition}

The following Proposition directly follows from Lemma~\ref{constantelemm} and It\^{o}'s formula. It characterizes martingales in terms of  the invariance of the finite variation part of the law of  processes by the associated lifted transformations of space depending on time:

\begin{proposition}Let $\nu\in \mathbb{S}$ and assume that $u : (t,x)\in [0,1] \times \mathbb{R}^d \to u_t(x)\in \mathbb{R}^d$ is such that the mapping $$h:= I_{\mathbb{R}^d}+  u$$ satisfy the hypothesis of Lemma~\ref{constantelemm}. Denote by $\Gamma$ the lift of $h$ (see  Definition~\ref{lift}) on $\mathbb{S}$ at $\nu$.  Then the following assertions are equivalent :\begin{enumerate}[(i)]\item $(t,\omega) \to u(t,W_t(\omega))$ is a $(\mathcal{F}_t^\nu)-$martingale\item $\nu-a.s.$ $$\int_0^. \left(\partial_t u_t(W_t) + (v_t^\nu. \nabla) u_t(W_t) + \sum_{i,j}\frac{{\alpha^\nu_t}^{i,j}}{2} \partial^2_{i,j}u(t,W_t) \right) dt = 0$$\item We have $\nu-a.s.$  \begin{equation} \label{transvconst} \int_0^. v_t^{{\Gamma}_\star\nu}\circ \Gamma dt = \int_0^. v_t^\nu dt \end{equation}  \end{enumerate}\end{proposition}

\section{Differential calculus associated to the variation processes}
\label{section4}
\label{differdef}

For $\nu \in \mathcal{P}(W)$ and $k\in V^{0}_\nu$ set $$\tau_{k}:= I_W + k$$ where $I_W$ stands for the $\nu-$equivalence class of mappings $\nu-a.s.$ equal to the identity on $W$. For $\epsilon \in \mathbb{R}$ set  $$\nu^k_\epsilon := {\tau_{\epsilon k}}_\star \nu = (I_W+ \epsilon 
k)_\star\nu$$ thus defining a path $$\epsilon \in \mathbb{R}\to \nu^k_\epsilon \in \mathcal{P}(W)$$  As we have seen  (Lemma~\ref{caspartp9})  if $\nu^k_0\in \mathbb{S}$ then 
for all $\epsilon  \in \mathbb{R}$, $\nu^k_\epsilon \in \mathbb{S}$ so that we have a path
$$\epsilon \in \mathbb{R}\to \nu^k_\epsilon \in \mathbb{S}$$  
on $\mathbb{S}$. Moreover  if $b^\epsilon$ 
denotes the finite variation part of $\nu^k_\epsilon\in \mathbb{S}$ (see Definition~\ref{defvelocity}), then  $b^\epsilon$ is simply related to $b^0$ (the finite variation part of $\nu$) 
by  $$b^\epsilon = b^0+ \epsilon k$$  Finally, since $\tau_{\epsilon k}$ is an isomorphism of filtered probability spaces, we can identify isometrically the spaces $L^p_a(\nu_\epsilon, H)$ along any such path $(\nu_\epsilon)$ as well as the vector spaces of variations processes (see Proposition~\ref{spapropisom}).  This motivates the following :

\begin{definition}
Given a mapping $$F : \nu \in \mathbb{S} \to F(\nu)\in \mathbb{R}\cup\{+\infty\}$$ and $\nu \in \mathbb{S}$ such that  $F(\nu)< \infty$, $F$  will be said to be $L^2_a(\nu, H_{0,0})$-differentiable (resp. $L^2_a(\nu, H)$-differentiable) at $\nu$ if for all $k\in V_\nu^{0,\infty}$ (resp. $k\in V_\nu^{\infty}$) $\frac{d}{d\epsilon} F(\nu_\epsilon^k)|_{\epsilon =0}$ exists
where for all $\epsilon \in \mathbb{R}$ and $k\in V_\nu^{0,\infty}$ (resp. $k\in V_\nu^{\infty}$), $$\nu_\epsilon^k := (\tau_{\epsilon k})_\star\nu := (I_W+ \epsilon k)_\star\nu ,$$
and if, in addition, there exists $\xi\in L^2_a(\nu, H)$  such that  for all $k\in V_\nu^{0,\infty}$ (resp. $k\in V_\nu^{\infty}$) 
\begin{equation} \label{mpmds} \frac{d}{d\epsilon} F(\nu_\epsilon^k)|_{\epsilon =0} = E_\nu\left[<\xi, k>_H\right] \end{equation} In this case, by Lemma~\ref{denselemma} (resp. Proposition~\ref{lemmatangent}) $\xi$ is the  unique element of $L^2_a(\nu, H_{0,0})$ (resp. in $ L^2_a(\nu, H)$), that we note $grad_\nu^0 F$ (resp. $grad_\nu F$), such that~(\ref{mpmds}) holds for any $k\in V_\nu^{0,\infty}$ (resp. in $V_\nu^\infty$), namely the orthogonal projection on $L^2_a(\nu,H_{0,0})$ of any (and then all) (resp. the unique element) $\xi$ satisfying~(\ref{mpmds}) for all $k\in V_\nu^{0,\infty}$ (resp. $k\in V_\nu^{\infty}$). We denote by $\delta F_\nu$ the linear continuous form defined by   \begin{equation} \label{dFdef} \delta F_\nu : h\in L^2_a(\nu, H_{0,0})  \to \delta F_\nu(h) := E_\nu[< grad_\nu^0 F, h>_H] \in \mathbb{R} \end{equation} (resp. by the same symbol $\delta F_\nu : h\in L^2_a(\nu, H)  \to \delta F_\nu(h) := E_\nu[< grad_\nu F, h>_H]\in \mathbb{R}$)
\end{definition}

\begin{remarkk}
 If $F$ is $L^2_a(\nu, H)-$differentiable at some $\nu\in \mathbb{S}$, then for all $k\in L^2_a(\nu,H_{0,0})$ we  have  $$E_\nu[< grad_\nu^0 F, k>_H] = E_\nu[< grad_\nu F, k>_H] $$ Indeed for $k\in L^2_a(\nu,H_{0,0})$ take $(k^n)\subset V_{\nu}^{0,\infty}$ a sequence which converges to $h$. We then have, for all $n\in \mathbb{N}$,
$$ E_\nu[< grad_\nu F, k^n>_H]  = \frac{d}{d\epsilon}F((I_W+ \epsilon k_n)_\star\nu)|_{\epsilon = 0} =  E_\nu[< grad_\nu^0 F, k^n>_H] $$ so that by continuity
$$ E_\nu[< grad_\nu F, k>_H] = \lim_{n\to \infty}E_\nu[< grad_\nu F, k^n>_H]  = \lim_{n\to \infty}E_\nu[< grad_\nu^0 F, k^n>_H]  = E_\nu[< grad_\nu^0 F, k>_H] $$
In particular  $ grad_\nu^0F$ is the orthogonal projection of  $ grad_\nu F$ on $L^2_a(\nu,H_{0,0})$ and it is meaningful to denote by the same symbol $\delta F_\nu$ the two linear forms. 
 \end{remarkk}

\section{The stochastic extension of Hamilton's least action principle}
\label{section5}
\subsection{Regular Lagrangians and their actions}

\begin{definition}
\label{defregular0}
A Borel measurable mapping $$\mathcal{L} : (t,x,v,a) \in [0,1]\times \mathbb{R}^d\times \mathbb{R}^d\times \mathcal{M}_d(\mathbb{R}) \to \mathcal{L}_t(x,v,a)\in \mathbb{R} \cup\{+\infty\}$$ will 
be called a \textbf{Lagrangian}. We denote its domain by $$Dom(\mathcal{L}) := \{ (t,x,v,a) : \mathcal{L}_t(x,v,a) < \infty\}$$ And we define the \textbf{action} of $\mathcal{L}$ on $\mathbb{S}$ to be the mapping $$\mathcal
{S} : \nu \in \mathbb{S} \to \mathcal{S}(\nu) \in \mathbb{R} \cup\{+\infty\}$$ defined by $$\mathcal{S}(\nu) := E_\nu\left[ \int_0^1 \mathcal{L}_t(W_t, v_t^\nu,\alpha_t^\nu) dt \right] $$ if 
$E_\nu\left[ \int_0^1| \mathcal{L}_t(W_t, v_t^\nu,\alpha_t^\nu) |dt \right] <\infty$ and by $\mathcal{S}(\nu)=\infty$ otherwise. 
\end{definition}
The following conditions will be used to ensure the least action principle : 

\begin{definition}
\label{defregular}
A Lagrangian $\mathcal{L}$ (see Definition~\ref{defregular0}) will be said to be \textbf{regular} if it satisfies the following assumptions
\begin{enumerate}[(i)]
\item The domain of $\mathcal{L}$ is the whole space i.e.  $$Dom(\mathcal{L})= [0,1]\times \mathbb{R}^d\times \mathbb{R}^d\times \mathcal{M}_d(\mathbb{R})$$
\item For all $(t,x,v,a)\in Dom(\mathcal{L})$, the mapping $$\mathcal{L}(t,x,v,a) : (\widetilde{x}, \widetilde{v})\in \mathbb{R}^d\times \mathbb{R}^d  \to \mathcal{L}_t(x+\widetilde{x},v+ \widetilde{v},a) \in \mathbb{R}$$ is Fr\'{e}chet differentiable at $(0_{\mathbb{R}^d},0_{\mathbb{R}^d})$ and we denote by \begin{equation}  \label{DLdef} D\mathcal{L}_{t,x,v,a} : (\widetilde{x}, \widetilde{v}) \in  \mathbb{R}^d 
\times \mathbb{R}^d\to  <(\partial_x \mathcal{L}_t)(x,v,a), \widetilde{x}>_{\mathbb{R}^d}  + <(\partial_v \mathcal{L}_t)(x,v,a), \widetilde{v}>_{\mathbb{R}^d}\end{equation} its 
derivative which the linear operator defined by  \begin{equation} \label{mdlm} D\mathcal{L}_{t,x,v,a}[\widetilde{x}, \widetilde{v}] := \frac{d}{d\epsilon}\mathcal{L}_t(x+ \epsilon \widetilde{x} ,v + \epsilon \widetilde{v} ,a )|_{\epsilon =0}\end{equation} 
\item The mappings $(t,x,v,a)\in Dom(\mathcal{L}) \to \partial_x \mathcal{L}_t(x,v,a)\in \mathbb{R}^d$ and $(t,x,v,a)\in Dom(\mathcal{L})  \to \partial_v \mathcal{L}_t(x,v,a)\in \mathbb{R}^d$ are Borel measurable.
\end{enumerate}
\end{definition}

\subsection{The stochastic least action principle and the related Euler-Lagrange condition}

\begin{definition}
\label{stoel}
 Let  $\mathcal{L}: (t,x,v,a) \in [0,1] \times \mathbb{R}^d\times \mathbb{R}^d \times \mathcal{M}_d(\mathbb{R}) \to \mathcal{L}_t(x,v,a) \in \mathbb{R}$  be a regular Lagrangian (see Definition~\ref{defregular}). A probability  $\nu \in \mathbb{S}$  such that for all $T<1$,  $\nu-a.s.$,
 $$\int_0^T| \partial_x \mathcal{L}_s(W_s, v_s^\nu, \alpha_s^\nu) |_{\mathbb{R}^d} ds  <\infty$$ 
 is said to  satisfy the \textbf{stochastic Euler-Lagrange condition} if and only if 
 there exists a $\mathbb{R}^d-$valued \textit{c\`{a}dl\`{a}g} $(\mathcal{F}_t^\nu)-$martingale $(N^\nu_t)_{t\in[0,1)}$ on the probability space $(W,\mathcal{B}(W)^\nu,\nu)$  such that  $\nu\otimes \lambda-a.s.$
\begin{equation}
\label{sel}
\partial_v \mathcal{L}_t(W_t, v_t^\nu, \alpha_s^\nu) - \int_0^t \partial_x \mathcal{L}_s(W_s, v_s^\nu, \alpha_s^\nu) ds  := N^\nu_t
\end{equation}
where $(v_t^\nu)$ (resp. $(\alpha_t^\nu)$) are the derivatives of the characteristics of $\nu$ (see Definition~\ref{defvelocity}). 
\end{definition}

\begin{theorem}
\label{thmlap}
Let $\mathcal{L}$ be a regular Lagrangian whose associated action on $\mathbb{S}$ is noted $\mathcal{S}$ (see Definitions~\ref{defregular0} and~\ref{defregular}). Assume also the existence of a strictly positive continuous function $f :\mathbb{R}^d \to \mathbb{R}^+$ and $p_1,p_2\geq 2$ such that
 \begin{equation} \label{H2lap} \limsup_{|\epsilon| \downarrow 0} \sup_{(t, x, v,a, \widetilde{x}, \widetilde{v}) \in Dom(\mathcal{L}) \times \mathbb{R}^d\times \mathbb{R}^d} \left( \frac{\left| \mathcal{L}_t(x+\epsilon \widetilde{x}, v+ \epsilon\widetilde{v}, a) -  \mathcal{L}_t(x,v,a) - \epsilon D\mathcal{L}_{t,x,v,a}[\widetilde{x}, \widetilde{v}]  \right|}{\epsilon f(\widetilde{x})\left(1+ |\widetilde{v}|_{\mathbb{R}^d}^2 + G(t,x,v,a)  \right)}\right)=0 \end{equation} holds, where $$G(t,x,v,a):= |\mathcal{L}_t(x,v,a)| + |\partial_x \mathcal{L}(x,v,a)|_{\mathbb{R}^d}^{p_1}  + |\partial_v \mathcal{L}(x,v,a)|_{\mathbb{R}^d }^{p_2} $$ Then for any $\nu \in \mathbb{S}$ satisfying 
 \begin{equation} \label{H0lap} \mathcal{S}(\nu) + E_\nu\left[\int_0^1 |\partial_x \mathcal{L}(W_s, v_s^\nu, \alpha_s^\nu)|^{p_1}_{\mathbb{R}^d}dt \right] + E_\nu\left[\int_0^1 |\partial_v \mathcal{L}(W_s, v_s^\nu, \alpha_s^\nu)|^{p_2}_{\mathbb{R}^d}dt \right] <\infty \end{equation} we have that $\mathcal{S}$ is $L^2_a(\nu,H_{0,0})-$differentiable (see subsection~\ref{differdef}) at $\nu$. Moreover $\nu$ satisfies the  stochastic Euler$-$Lagrange condition (see Definition~\ref{stoel}) if and only if  $$\delta \mathcal{S}_\nu(h) =0$$ for all $h\in L^2_a(\nu, H_{0,0})$ i.e. if and only if for all $h\in V^{0,\infty}_\nu$, $$\frac{d\mathcal{S}({\tau_{\epsilon h}}_\star\nu)}{d\epsilon}|_{\epsilon=0} =0 $$where for $h\in V^{0,\infty}_\nu$ $$\tau_h := I_W + h$$  \end{theorem}
\nproof
By~(\ref{H0lap}) since $p_1,p_2 \geq 2$ we have \begin{equation}\label{lemmahold} E_\nu\left[\int_0^1 |\partial_v \mathcal{L}_t(W_t, v_t^\nu, \alpha_t^\nu) |^2 dt \right]  + E_\nu\left[\int_0^1 |\partial_x \mathcal{L}_t(W_t, v_t^\nu, \alpha_t^\nu) |^2 dt \right]  <\infty \end{equation} Define $$\dot{\xi}_t:=  \partial_v \mathcal{L}_t(W_t, v_t^\nu, \alpha_t^\nu)- \int_0^t \partial_x \mathcal{L}_s(W_s, v_s^\nu, \alpha_s^\nu)ds$$
and $$\xi:= \int_0^. \dot{\xi}_t dt.$$  From ~(\ref{lemmahold}) and Jensen's inequality, $\xi \in L^2_a(\nu,H)$. Take $h\in V_\nu^{0,\infty}$ and set
$$A_\epsilon := \left| \frac{\mathcal{S}({\tau_{\epsilon h}}_\star \nu) - \mathcal{S}(\nu)}{\epsilon} - E_\nu[<h,\xi>_H] \right|$$ where ${\tau_{\epsilon h}}_\star \nu$ is the image of the probability $\nu$ by the mapping $\tau_{\epsilon h} := I_W + \epsilon h$.  We want to show that $A_\epsilon$ converges to $0$. By Lemma~\ref{caspartp9} we have 
\begin{eqnarray*}  S({\tau_{\epsilon h}}_\star \nu)  & = & E_{{\tau_{\epsilon h}}_\star \nu}\left[ \int_0^1 \mathcal{L}_s(W_s, v^{{\tau_{\epsilon h}}_\star \nu}_s, \alpha^{{\tau_{\epsilon h}}_\star \nu}_s) ds  \right]  \\ & = & E_{ \nu}\left[ \int_0^1 \mathcal{L}_s(W_s+ \epsilon h_s, v^{{\tau_{\epsilon h}}_\star \nu}_s\circ \tau_{\epsilon h}, \alpha_s^{{\tau_{\epsilon h}}_\star \nu} \circ { \tau_{\epsilon h}}) ds \right]  \\ & =&  E_{ \nu}\left[ \int_0^1 \mathcal{L}_s(W_s + \epsilon h_s, v_s^\nu + \epsilon \dot{h}_s, \alpha^\nu_s) ds \right] 
\end{eqnarray*}
so that we first obtain, for $\epsilon \in \mathbb{R}$,
$$A_\epsilon = \left| E_\nu\left[\int_0^1 \left(\frac{\mathcal{L}_s(W_s + \epsilon h_s, v_s^\nu + \epsilon \dot{h}_s, \alpha^\nu_s) -\mathcal{L}_s(W_s, v_s^\nu, \alpha^\nu_s)}{\epsilon} -  <\dot{h}_s,\dot{\xi}_s>\right)ds\right]\right|$$
On the other hand, since $\nu-a.s.$ $h_0= h_1=0$, an integration yields $$A_\epsilon = \left| E_\nu\left[\int_0^1 \left( \frac{\mathcal{L}_s(W_s + \epsilon h_s, v_s^\nu + \epsilon \dot{h}_s, \alpha^\nu_s) -\mathcal{L}_s(W_s, v_s^\nu, \alpha^\nu_s)}{\epsilon} - D\mathcal{L}_t(W_s, v_s^\nu, \alpha_s^\nu)[h_s, \dot{h}_s]   \right)ds \right]\right|$$ where $(t,x,v,a)\to  D\mathcal{L}_t(x,v,a)$ is given by~(\ref{DLdef}).  By the hypothesis $$0 \leq E_\nu\left[\int_0^1 G(t,W_t,v_t^\nu,\alpha_t^\nu)dt \right]<\infty$$ For $\epsilon_0 >0$ by~(\ref{H2lap})  there exists  $\alpha >0$ such that, for all $\epsilon\in \mathbb{R} / \{0\}$ with $|\epsilon |< \alpha$, the following inequality holds everywhere
$$\left|\frac{\mathcal{L}_t(x+ \epsilon \widetilde{x}, v+ \epsilon \widetilde{v}, a) -\mathcal{L}_t(x,v,a)}{\epsilon} - D\mathcal{L}_t(x, v a )[\widetilde{x}, \widetilde{v}]\right| \leq \frac{\epsilon_0}{B} f(\widetilde{x})(1+ |\widetilde{v}|^2_{\mathbb{R}^d}+ G(t,x,v,a) )$$  where $$B:=(\sup_{z\in B(0_{\mathbb{R}^d}, R) } f(z)) \left(1+ |h|_{L^2(\nu,H)}^2 +E_\nu\left[\int_0^1 G(t,W_t,v_t^\nu,\alpha_t^\nu) dt \right]\right)$$ and    $$R:= |h|_{L^{\infty}(\mu,W)}+1 <\infty$$ Since $h\in V^{0,\infty}_\nu$ we have $\nu-a.s. $ for all $t\in[0,1]$, $|h_t(\omega)|_{\mathbb{R}^d}\leq |h|_{L^\infty(\nu,W)}$. Thus by the continuity of $f$, for any $\epsilon$ such that $|\epsilon | <\alpha$  we obtain \begin{eqnarray*}  A_\epsilon & \leq & \frac{\epsilon_0}{B} E_\nu\left[\int_0^1 f(h_s)(1+ |\dot{h}_s|_{\mathbb{R}^d}^2 + G(W_s, v_s^\nu, \alpha_s^\nu) )ds \right]
\\ & \leq & \frac{\epsilon_0}{B}   E_\nu\left[\int_0^1 \left(\sup_{z\in B(0_{\mathbb{R}^d}, R)}f(z)\right)(1+ |\dot{h}_s|_{\mathbb{R}^d}^2 + G(W_s, v_s^\nu, \alpha_s^\nu)  )ds \right]
\\& \leq & \frac{\epsilon_0}{B}   \left(\sup_{z\in B(0_{\mathbb{R}^d}, R)}f(z)\right)  E_\nu\left[\int_0^1 (1+ |\dot{h}_s|_{\mathbb{R}^d}^2 +G(W_s, v_s^\nu, \alpha_s^\nu) )  ds \right]
\\ & \leq & \epsilon_0 
\end{eqnarray*} Hence we have $$\limsup_{|\epsilon|\downarrow 0} A_\epsilon \leq \epsilon_0$$ Since this inequality holds for any $\epsilon_0>0$, we conclude that $\limsup_{|\epsilon|\downarrow 0} A_\epsilon =0$. Using the definition of $A_\epsilon$ we get that for all $h\in V_\nu^{0,\infty}$ $\frac{d}{d\epsilon}\mathcal{S}{(\tau_{\epsilon h}}_\star \nu) |_{\epsilon =0}$ exists and
\begin{equation} \label{pmastar}\delta \mathcal{S}_\nu[h] = \frac{d}{d\epsilon}\mathcal{S}({\tau_{\epsilon h}}_\star \nu) |_{\epsilon =0} = E_\nu [<\xi, h>_H] \end{equation} In particular, $\mathcal{S}$ is $L^2_a(\nu,H_{0,0})-$differentiable (see subsection~\ref{differdef}). Note that by Proposition~\ref{emery} and by definition of $\xi$,  $\nu$ satisfies the stochastic Euler$-$Lagrange condition if and only if $\xi$ is orthogonal to $L^2_a(\nu,H_{0,0})$. On the other hand, by~(\ref{pmastar})  $\nu$ satisfies the least action principle for $\mathcal{S}$ if and only if  for all $h\in V_\nu^{0,\infty}$ $$E_\nu\left[< h, \xi>_H\right] =0$$ while by Lemma~\ref{denselemma} this latter condition is also satisfied if and only if $\xi$ is orthogonal to $L^2_a(\nu,H_{0,0})$. This achieves the proof.
 \nqed
\begin{remarkk}
\begin{itemize}
\item Whenever  $\nu$ satisfies the stochastic Euler-Lagrange condition it also satisfies the following averaged stochastic Euler-Lagrange condition~(\ref{aver})  
\begin{equation}
\label{aver}
\frac{d}{dt}E_\nu\left[\partial_v \mathcal{L}_t(W_t, v_t^\nu, \alpha_t^\nu)\right] = E_\nu\left[ \partial_x \mathcal{L}_s(W_t, v_t^\nu, \alpha_t^\nu)\right]
\end{equation}
Moreover  the left hand side of~(\ref{aver}) is well defined and it is trivial to check that~(\ref{aver}) holds if and only if  for any $h\in H_{0,0}$ we have 
\begin{equation} \label{determ} \delta \mathcal{S}_\nu[h]=0\end{equation}
\item We refer to \cite{RLZ} for an Hamiltonian point of view on the stochastic Euler$-$Lagrange condition.
\item Assume that $\mathcal{L}$ is a regular Lagrangian which does not depends on $(\alpha_s^\nu)$. For $x\in \mathbb{R}^d$ and $u\in H$ let $\gamma\in W$ be defined by $\gamma_t:= x + u_t$ for all $t\in
[0,1]$, and denote by $\delta^{Dirac}_\gamma$ the Dirac measure concentrated on $\gamma$, namely  $\delta^{Dirac}_\gamma(A)= I_A(\gamma)$ for $A\in \mathcal{B}(W)$. Then $(\mathcal{F}_t^\nu)$ is the filtration constant equal 
to the set of the subsets of $W$, so that the martingales can be identified with the constants, and $V_\nu^{0,\infty} = V_\nu^0 = L^2_a(\nu,H_{0,0})\simeq H_{0,0}$. Moreover in this case for $h\in V_\nu^0$, $\epsilon \in \mathbb
{R}$,  we have ${\tau_{\epsilon h}}_\star \nu = \delta_{\gamma +\epsilon h}$ so that  ${\tau_{\epsilon h}}_\star \nu \in \mathbb{S}$ has a martingale part equal to $0$  and $ {\tau_{\epsilon h}}_\star \nu$ a.s. for all $t\in[0,1]$,  $$W_t= 
x + \gamma_t + \epsilon h_t$$ In particular $\dot{W}_t$ exists a.s. and $ {\tau_{\epsilon h}}_\star \nu$ a.s. $v_t^{{\tau_{\epsilon h }}_\star\nu }(\omega)=\dot{W}_t= \dot{\gamma}_t + \epsilon \dot{h}_t$. Hence we obtain  $\mathcal
{S}({\tau_{\epsilon h}}_\star \nu) = \int_0^1 \mathcal{L}(\gamma_s + \epsilon h_s, \dot{\gamma}_s + \epsilon \dot{h}_s)ds$ and the stochastic least action principle reads  $$\frac{d}{d\epsilon }|_{\epsilon= 0}\int_0^1 \mathcal{L}(\gamma_s + \epsilon h_s, \dot{\gamma}_s + \epsilon \dot{h}_s)=0$$ for all $h\in H_{0,0}$, while the stochastic 
Euler-Lagrange condition holds if and only if there exists  $c\in \mathbb{R}$ such that $\lambda-a.s.$ $$\partial_v \mathcal{L}(\gamma_t, \dot{\gamma}_t)- \int_0^t\partial_x \mathcal{L}(\gamma_s, \dot{\gamma}_s) ds = c$$ or if 
and only if $t\to \partial_v \mathcal{L}(t,\gamma_t, \dot{\gamma}_t)$ is differentiable a.e. and $\lambda-a.e..$  \begin{equation} \label{elclass} \frac{d}{dt}\partial_v \mathcal{L}(\gamma_t, \dot{\gamma}_t) = \partial_x \mathcal{L}(\gamma_t, \dot{\gamma}_t) \end{equation}
Moreover, if we assume also that $(x,v) \in \mathbb{R}^d\times \mathbb{R}^d \to  \partial_x \mathcal{L}(x, v)\in \mathbb{R}^d$ is continuous (resp. $(x,v) \in \mathbb{R}^d\times \mathbb{R}^d \to  \partial_v \mathcal{L}(x, v)\in \mathbb{R}^d$is $C^1$), and that $t\in[0,1] \to \gamma_t \in \mathbb{R}^d$ is $C^2$, then~(\ref{elclass}) holds for all $t\in[0,1]$. In this case we recover the least action principle of classical mechanics.
\end{itemize}
\end{remarkk}

\section{Invariances and Noether's theorem}
\label{section6}

\begin{definition}
\label{invariant}
Let $h : [0,1]\times \mathbb{R}^d\to \mathbb{R}^d$ be a mapping which satisfies the hypothesis of Lemma~\ref{constantelemm}. We say that $h$ is  an $\mathbb{S}-$\textbf{invariant transformation} for $\mathcal{L}$ if for all $\nu\in \mathbb{S}$, $\nu-a.s $ all $\omega\in W $ we have,  $\lambda-a.s.$,  \begin{equation} \label{noethearnold} \mathcal{L}_t(\Gamma_t(\omega), v_t^{{\Gamma}_\star \nu}\circ \Gamma
(\omega), \alpha_t^{{\Gamma}_\star \nu}\circ \Gamma(\omega))  =   \mathcal{L}_t(W_t(\omega), v_t^{\nu}(\omega), \alpha_t^{\nu}(\omega))\end{equation}
where $\Gamma$ is the lift of $h$ on $\mathbb{S}$ at $\nu$ (see  Definition~\ref{lift}). Moreover we say that  a family $(h^\epsilon)_{\epsilon \in \mathbb{R}}$ of $\mathbb{S}-
$invariant transformation for $\mathcal{L}$ is a \textbf{differentiable family of $\mathbb{S}-$invariant transformations} for $\mathcal{L}$, if     $(t,x,\epsilon)\in [0,1]\times \mathbb{R}^d\times \mathbb{R} \to h^\epsilon(t,x)\in \mathbb{R}^d$ is $C^1$ in $\epsilon$ and $h^0(t,x)= x$ for all $(t,x)\in [0,1]\times \mathbb{R}^d$.  
\end{definition}  

We recall that for two real valued \textit{c\`{a}dl\`{a}g} semi-martingales $X$ and $Y$, their quadratic co-variation is the process ($[X,Y]_t$) defined by $$[X,Y]_. = X_.Y_.- \int_0^. X_{s-}dY_s -\int_0^. Y_{s-}dX_s$$ see \cite{DM2} for more.

\begin{theorem}
\label{noether}
Let $\mathcal{L}$ be a regular Lagrangian which is $C^1$ and assume that $\nu \in \mathbb{S}$ satisfies the stochastic Euler-Lagrange conditions (see Definition~\ref{stoel}) for $\mathcal{L}$.  Assume also that $(h^\epsilon)_{\epsilon \in \mathbb{R}}$ is a differentiable family of $\mathbb{S}-$invariant transformations for $\mathcal{L}$. Let $(\mathcal{I}_t)_{t\in[0,1]}$ be any optional process on $(W,\mathcal{F}^\nu_.,\nu)$ such that $\lambda\otimes \nu-a.s.$
\begin{equation} \label{noetherconst} \mathcal{I}_t:= <\frac{d}{d\epsilon}|_{\epsilon=0}h^\epsilon_t(W_t), p_t^\nu >_{\mathbb{R}^d} - \sum_i \left[ \frac{d}{d\epsilon}|_{\epsilon=0}{h^\epsilon_.(W_.)}^i, {p_.^\nu}^i \right]_t + \int_0^t  \theta_s ds \end{equation} where $[. ,. ]$ denotes the quadratic co-variation of semi-martingales, $(p_t^\nu)$ denotes a \textit{c\`{a}dl\`{a}g} modification of the process $(\partial_v \mathcal{L}_t(W_t,v_t^\nu, \alpha_t^\nu))$,
and
\begin{equation} \label{thetas} \theta_s :=  \sum_{i,j}{{\kappa}_s^{i,j}} \frac{\partial \mathcal{L}}{\partial \alpha_{i,j}}(W_s,v_s^\nu, \alpha_s^\nu) \end{equation} and where  $
(\kappa_s(\omega))$ is the $\mathcal{M}_d(\mathbb{R})-$valued process defined by  \begin{equation} \label{kappadl} {\kappa}_s(\omega):=  {\alpha^\nu_s}.
\left(\left(\nabla \frac{d}{d\epsilon}  h^\epsilon|_{\epsilon=0}\right)(s,W_s) \right)^\dagger +  \left(\left(\nabla \frac{d}{d\epsilon}  h^\epsilon|_{\epsilon=0}\right)(s,W_s)\right)  {\alpha^
\nu_s} \end{equation}   Then $\left(\mathcal{I}_t \right)_{t\in[0,1)}$ is a $(\mathcal{F}_t^\nu)-$local martingale on the probability space $(W, \mathcal{B}(W)^\nu, \nu)$.
\end{theorem}
\nproof
For all $\epsilon \in \mathbb{R}$ we denote by $\Gamma^\epsilon$ the lift of $h^\epsilon$ on $\mathbb{S}$ at $\nu$ (see Definition~\ref{lift}). We set $\nu_\epsilon:= {\Gamma^{\epsilon}}_\star\nu$ and $$\widetilde{u}(t,x):= \frac{d}{d\epsilon}h^\epsilon_t(x)|_{\epsilon=0}$$ For $T<1$, by Lemma~\ref{constantelemm},
\begin{equation} 
\label{lplo}
\int_0^T \mathcal{L}(\Gamma^\epsilon_t, v_t^{\nu_\epsilon}\circ\Gamma^\epsilon, \alpha_t^{\nu_\epsilon}\circ \Gamma^\epsilon) dt = \int_0^T \mathcal{L}(h^\epsilon(t,W_t),  m^\epsilon_t,   a_t^\epsilon) dt 
\end{equation}
where $$m^\epsilon_t= \partial_t {h}^\epsilon(t,W_t) + (v_t^\nu. \nabla) {h}^\epsilon(t,W_t) + \sum_{i,j}\frac{{\alpha^\nu_t}^{i,j}}{2} \partial^2_{i,j}{h}^\epsilon(t,W_t) $$
and  $$a_t^\epsilon = \alpha_t^\nu .(\nabla {h^\epsilon})^\dagger(t,W_t) + \nabla {h^\epsilon}(t,W_t).\alpha_t^\nu $$
By differentiating~(\ref{lplo}), condition~(\ref{noethearnold}) yields $\nu-a.s.$,
\begin{equation}
\label{ltpremres}
0= \int_0^T \left( <\widetilde{u}(t,W_t),\partial_x \mathcal{L}(W_t, v_t^\nu,\alpha_t^\nu) > + <Q_t,\partial_v \mathcal{L}(W_t, v_t^\nu,\alpha_t^\nu) > + \theta_t  \right) dt
\end{equation} where \begin{equation} \label{Qt} Q_t := \partial_t \widetilde{u} (t,W_t) + (v_t^\nu. \nabla) \widetilde{u} (t,W_t) +  \sum_{i,j}\frac{{\alpha^\nu_t}^{i,j}}{2} \partial^2_{i,j}{\widetilde{u}}(t,W_t) 
\end{equation}  and where $(\theta_s)$ is given by~(\ref{thetas}) and $\kappa$ by~(\ref{kappadl}).  Since $\nu$ satisfies the stochastic Euler$-$Lagrange condition there exists a $(\mathcal{F}_t^\nu)$ \textit{c\`{a}dl\`{a}g} martingale $(N_t^\nu)$ such that $\nu\otimes \lambda-a.s.$ \begin{equation} \label{elenoe} \partial_v \mathcal{L}_t(W_t,v_t^\nu, \alpha_t^\nu) =  N_t^\nu  +\int_0^t \partial_x \mathcal{L}_s(W_s,v_s^\nu, \alpha_s^\nu) ds\end{equation} Denote by $(p_t^\nu)$ the process defined by the right hand term of~(\ref{elenoe}). We have  \begin{equation} \label{internoet2} \int_0^T  <\widetilde{u}(t,W_t), \partial_x \mathcal{L}(W_t, v_t^\nu,\alpha_t^\nu) > dt  =   \int_0^T <\widetilde{u}(t,W_t),dp_t^\nu>   - \int_0^T <\widetilde{u}(t,W_t),dN_t^\nu>  \end{equation} We now compute the first term of the right hand side. Denoting by $M^\nu$ the martingale part of $\nu$ by It\^{o}'s formula we obtain a.e. $$\widetilde{u}^i(t,W_t)= {M_t^{u,i}} + {A_t}^i $$ where $$M_t^{u,i}= \widetilde{u}^i(0,W_0) + \int_0^t \sum_j (\partial_j \widetilde{u}^i)(s,W_s) dM_s^{\nu,j}$$  and where
$$A_t^i := \int_0^t Q_s^i ds$$  with $(Q_s)$ given by~(\ref{Qt}). Since $(M_t^u)$ is a continuous local martingale and $(A_t)$ is continuous and of finite variation, by It\^{o}'s formula  (by (18.1) and (19.2) VIII p.343 of \cite{DM2}) 
we obtain \begin{multline}
\label{internoet3}
\int_0^T <\widetilde{u}(t, W_t),dp_t^\nu> = <\widetilde{u}(T,W_T),p_T^\nu>_{\mathbb{R}^d}   - \int_0^T <p_{t_{-}}^\nu,dM_t^u> ... \\...  -\int_0^T <p_t^\nu, Q_t> dt - \sum_i [M^{u,i}, {p^\nu}^i]_T 
\end{multline}Putting~(\ref{internoet3}) into~(\ref{internoet2}) we derive\begin{multline}
\label{milm}
 \int_0^T <\widetilde{u}(t,W_t), \partial_x \mathcal{L}(W_t, v_t^\nu,\alpha_t^\nu) > dt  =  <\widetilde{u}_T(W_T),p_T^\nu>_{\mathbb{R}^d} ...  \\ ...    -\int_0^T <p_t^\nu, Q_t> dt - \sum_i [M^{u,i}, {p^\nu}^i]_T - \widetilde{M}_T
\end{multline} 
where $$ \widetilde{M}_t =   \int_0^T <{p_{t-}}^\nu,dM_t^u>  + \int_0^T <\widetilde{u}(t,W_t), d{N_t}^\nu>   $$ By putting~(\ref{milm}) into~(\ref{ltpremres}) we obtain \begin{equation} \label{finalnoeth} \mathcal{I}_T= \widetilde{M}_T\end{equation} On the other hand  by construction $(\widetilde{M}_t)$ is a \textit{c\`{a}dl\`{a}g} $(\mathcal{F}_t^\nu)$-local martingale so that the result follows by~(\ref{finalnoeth}).
\nqed

\begin{remarkk}
This theorem must be compared with the original theorem such as it is formulated  p.88 of  \cite{ARNOLD6} together with the following remark. Take $c\in W$ to be such that $t\in[0,1] \to c_t\in \mathbb{R}^d$ is smooth, $\nu:= \delta_c$ the associated Dirac measure and the probability space $(W,\mathcal{B}(W)^\nu, \nu)$. Consider a particular $h: (t,x) \in [0,1] \times \mathbb{R}^d \to h(x) \in \mathbb{R}^d$ not depending on time  and satisfying the hypothesis of Lemma~\ref{constantelemm}. and let $\Gamma$ denote the lift of $h$ at $\nu$. Then for this transformation, we have ${\Gamma}_\star\nu = \delta_{\widetilde{c}}$ where $\widetilde{c}_t = h(c_t)$ for all $t$, and $\lambda \otimes \nu-a.s.$ $$v_t^{{\Gamma}_\star\nu}\circ \Gamma= \frac{dh(W_t)}{dt}$$ The r.h.s. of this latter equation is nothing but the image of $\dot{W}_t$ by $h$ (noted $ h_\star (\dot{W}_t)$ in geometry).
\end{remarkk}

\section{Application to stochastic control}
\label{section7}

\begin{definition}
\label{stable}
A subset $\mathcal{N}$ of $\mathbb{S}$ is said to be $V^0-$stable if for any $\nu\in \mathcal{N}$ any $\epsilon \in \mathbb{R}$ and $h\in V_\nu^{0}$
 $${\tau_{\epsilon h}}_\star \nu \in \mathcal{N}$$ where $$\tau_{\epsilon h} := I_W+\epsilon h$$
  \end{definition}

Note that in particular (see Lemma~\ref{caspartp9}), $\mathbb{S}$ is $V^0-$stable. 

\begin{theorem}
\label{corollairecontrole0}
Let $\mathcal{N}$ be a $V^0-$stable subset of $\mathbb{S}$ (see Definition~\ref{stable}). Consider $\mathcal{L}$ a non negative  regular Lagrangian with associated action $\mathcal{S}$. Assume also that there exists a strictly positive continuous function $f: \mathbb{R}^d \to \mathbb{R}^+$ and $p_1,p_2\geq 2$ such that~(\ref{H2lap}) holds for $\mathcal{L}$  and  \begin{equation}  \label{H1lap} \sup_{(t,x,v,a) \in Dom(\mathcal{L})} \left( \frac{|\partial_x \mathcal{L}_t(x,v,a)|_
{\mathbb{R}^d}^{p_1} +  |\partial_v \mathcal{L}_t(x,v,a)|_{\mathbb{R}^d}^{p_2} )}{1 + \mathcal{L}_t(x,v,a)}\right)<\infty \end{equation}  Consider the minimization problem \begin{equation} \label{gammavar0} I_F:= \inf\left( \{ \mathcal{S}(\nu):  \nu\in \mathcal{N}\} \right) \end{equation} 
and assume that $I_F<\infty$. Then for any $\eta \in \mathcal{N}$ which attains the infimum of~(\ref{gammavar0}),  $\eta$ satisfies the stochastic Euler-Lagrange condition~(\ref{sel}) for $\mathcal{L}$. Moreover if $\mathcal{L}$ is $C^1$ and $(h^\epsilon)$ is a differentiable family of $\mathbb{S}-$invariant transformations for $\mathcal{L}$ (see Definition~\ref{invariant}), the process defined on $(W, \mathcal{B}(W)^\eta, \eta)$ by~(\ref{noetherconst}) is a $(\mathcal{F}_t^\eta)-$local martingale. Moreover the same statements hold if we change the $inf$ by $sup$, assuming it is attained.
\end{theorem}
\nproof
First notice that by~(\ref{H1lap}), since  $\mathcal{S}(\eta) <\infty$, ~(\ref{H0lap})  is satisfied with $p_1,p_2$ as~(\ref{H2lap}) is assumed. Thus Theorem~\ref{thmlap} applies to $\mathcal{L}$ at $\eta$ and we first obtain that   $\mathcal{S}$ is $L^2_a(\eta, H_{0,0})-$differentiable. Notice that, for any $h\in V_{\eta}^{0,\infty}$, by Definition~\ref{stable} we have $${\tau_h}_\star \eta \in  \mathcal{N}$$ Therefore since $\eta$ (it is a minimum of the action, it satisfies $$\delta \mathcal{S}_\eta[h] :=\frac{d}{d\epsilon} \mathcal{S}({\tau_{\epsilon h}}_\star\eta) |_{\epsilon=0} = 0$$ for all $h\in V_\eta^{0,\infty}$, so that by definition of $\delta \mathcal{S}_\eta$ we get $$\delta \mathcal{S}_\eta =0$$ Thus, the result directly follows by Theorem~\ref{thmlap} and Theorem~\ref{noether}. If we change the $inf$ by $sup$ the proof is similar.
\nqed

For $\gamma \in \mathcal{P}_{\mathbb{R}^d\times \mathbb{R}^d}$ a Borel probability on $\mathbb{R}^d\times \mathbb{R}^d$ (resp. $\nu_0,\nu_1$ two Borel probabilities on $\mathbb{R}^d$) we set $$\mathbb{S}_\gamma:= \{\nu \in \mathbb{S} : (W_0\times W_1)_\star \nu = \gamma\}$$  where $(W_0\times W_1): \omega\in W \to (W_0(\omega), W_1(\omega)) \in \mathbb{R}^d\times \mathbb{R}^d$  i.e. $\mathbb{S}_\gamma$ is the set of the $\nu \in \mathbb{S}$ with a fixed joint law $\gamma$ for $(W_0,W_1)$ Consider also $$\mathbb{S}_{\nu_0,\nu_1}:= \{\nu \in \mathbb{S}: {W_0}_\star \nu = \nu_0, {W_1}_\star \nu = \nu_1 \}$$ i.e. the set of $\nu\in \mathbb{S}$ with initial marginal $\nu_0$ and final marginal $\nu_1$. We also denote the set of the $\nu\in \mathbb{S}$ whose martingale part is a Brownian motion by $\mathbb{S}_B$ i.e.  $$\mathbb{S}_B :=\{ \nu \in \mathbb{S} : M^\nu_\star\nu = \mu^0\}$$ where $\mu^0$ is the standard Wiener measure (the law of the Brownian motion, with $\mu^0-a.s.$ $W_0=0$).  We refer to \cite{TOUZI} an to the proofs therein for sufficient condition for the existence to the following minimizer.  By assuming their existence we obtain :

 \begin{corollary}
 \label{corollarycontrole}
 Consider a probability $\gamma \in \mathcal{P}_{\mathbb{R}^d\times \mathbb{R}^d}$ (resp. $\nu_0,\nu_1 \in \mathcal{P}_{\mathbb{R}^d}$) and any non negative Lagrangian $\mathcal{L}$ which is regular and whose action is denoted by $\mathcal{S}$ (see Definition~\ref{defregular}). Assume also that  hypothesis~(\ref{H1lap}) and~(\ref{H2lap}) are satisfied for some $p_1,p_2\geq 2$  and some strictly positive continuous mapping $f: \mathbb{R}^d\to \mathbb{R}^+$.  Consider the variational problems    \begin{equation}  \label{gammavar} I_{\nu_0,\nu_1}:= \inf\left( \mathcal{S}(\nu):  \nu\in \mathbb{S}_{\nu_0,\nu_1} \right)  \end{equation}    \begin{equation}  \label{gammavar2}  I_\gamma :=  \inf\left( \mathcal{S}(\nu):  \nu\in \mathbb{S}_{\gamma} \right)\end{equation}  \begin{equation}  \label{gammavarB} I^B_{\nu_0,\nu_1}:= \inf\left( \mathcal{S}(\nu):  \nu\in \mathbb{S}_B \cap \mathbb{S}_{\nu_0,\nu_1} \right)  \end{equation}   and    \begin{equation}  \label{gammavar2B}  I^B_\gamma :=  \inf\left( \mathcal{S}(\nu):  \nu\in \mathbb{S}_B \cap \mathbb{S}_{\gamma} \right) \end{equation}   By assuming that $I_{\nu_0,\nu_1}$ (resp. $I_\gamma$, resp. $I_{\nu_0,\nu_1}^B$, resp.  $I_\gamma^B$) is finite (i.e. $<\infty$),  any $\nu \in \mathbb{S}_{\nu_0,\nu_1}$  (resp. $\mathbb{S}_\gamma$, resp. $\mathbb{S}_{\nu_0,\nu_1}\cap\mathbb{S}_B$, resp. $\mathbb{S}_\gamma\cap \mathbb{S}_B$) which attains  the infimum of~(\ref{gammavar}) (resp. of~(\ref{gammavar2}), resp. of (\ref{gammavarB}), resp. of (\ref{gammavar2B})), satisfies the stochastic Euler-Lagrange condition~(\ref{sel}) for $\mathcal{L}$. Moreover  if $\mathcal{L}$ is $C^1$ and $(h^\epsilon)$ is  a differentiable family of $\mathbb{S}-$invariant transformations for $\mathcal{L}$ (see Definition~\ref{invariant}) then the process defined on the probability space $(W, \mathcal{B}(W)^\nu, \nu)$ by~(\ref{noetherconst}) is a $(\mathcal{F}_t^\nu)-$local martingale. Moreover the same statements hold if we change the $inf$ by $sup$, assuming it is attained. 
  \end{corollary}

 \nproof
  Note that $ \mathbb{S}_{\gamma}$ is $V^0-$stable (see Definition~\ref{stable}). Indeed, for all $\nu \in \mathbb{S}_\gamma$ and  $h\in V^0_\nu$  we have, $\nu-a.s.$ $$h_0=h_1=0$$ so that $$({W_0\times W_1})_\star ({\tau_h}_\star \nu)= ((W_0+h_0)\times (W_1+ h_1))_\star \nu = (W_0\times W_1)_\star \nu = \gamma$$  Similarly $\mathbb{S}_{\nu_0,\nu_1}$  is $V^0-$stable.  Thus  Paul Levy's criterion and~(\ref{2R311}) of Lemma~\ref{caspartp9}  imply that $\mathbb{S}_\gamma\cap  \mathbb{S}_B$ and $\mathbb{S}_{\nu_0,\nu_1} \cap  \mathbb{S}_B$ are also $V^0-$stable. Therefore the results directly follow by Theorem~\ref{corollairecontrole0}.  For $sup$, the proof is similar.
\nqed

\section{The critical processes of the classical Lagrangians}
\label{section8}
\subsection{The classical Lagrangians}

\begin{definition} Given a function $V: [0,1]\times \mathbb{R}^d \to \mathbb{R}$ which is assumed to be measurable and $C^1$ in $x$, we define classical Lagrangians of the form  $$\mathcal{L}^V : (t,x,v,a)\in [0,1]\times \mathbb{R}^d\times \mathbb{R}^d \times \mathcal{M}_d(\mathbb{R}) \to \frac{|v|^2}{2} - V(x) \in \mathbb{R}$$  and denote by  $G_V(\mathbb{R}^d)$  the set of  $\nu \in \mathbb{S}$ which satisfy the stochastic Euler-Lagrange condition (see Definition~\ref{stoel}) for $\mathcal{L}^V$.
  \end{definition}

The following is the counterpart to the Galilean invariance for the free particle of classical mechanics in our stochastic framework :

\begin{proposition}A measure $\nu \in \mathbb{S}$ belongs to  $G_0(\mathbb{R}^d)$ if and only if for any $U\in \mathcal{I}_f^0(\nu)$  of the form $$U:=  I_W+  h$$ where $h\in\mathcal{M}_a(\nu, H) \cap V_\nu$ (see Proposition~\ref{emery} and Definition~\ref{defframe}) we have  $$U_\star \nu \in G_0(\mathbb{R}^d)$$\end{proposition}

\nproof For any $k:=\int_0^. \dot{k}_s ds \in L^2_a(U_\star \nu,H_{0,0})$, by Lemma~\ref{caspartp9} and Proposition~\ref{emery} we obtain \begin{eqnarray*} E_{{U}_\star \nu }\left[\int_0^1< v^{{U}_\star\nu}_s, \dot{k}_s> ds\right] & = & E_{\nu }\left[ \int_0^1< v^{{U}_\star\nu}_s\circ U, \dot{k}_s\circ U> \right] \\ & = & E_{\nu }\left[\int_0^1< v_s^\nu +\dot{h}_s  , \dot{k}_s\circ U> ds\right] \\ &  =  & E_{\nu }\left[\int_0^1< v_s^\nu, \dot{k}_s\circ U> ds\right]  + E_\nu\left[< h,k\circ U>_H\right]  \\ &  =  & E_{\nu }\left[\int_0^1< v_s^\nu, \dot{k}_s\circ U> ds\right]  \end{eqnarray*} On the other hand since $U$ is an isomorphism of filtered probability spaces, we have $k\in L^2_a(\nu,H_{0,0})$  iff there exists a $\widetilde{k}\in L^2_a(U_\star \nu,H_{0,0})$  such that $\nu-a.s.$ $k=\widetilde{k}\circ U$. Hence the result follows from Proposition~\ref{emery}. \nqed

Next proposition characterizes the measures in $G_V(\mathbb{R}^d)$.

\begin{proposition}
\label{sdegd}For $\nu \in G_V(\mathbb{R}^d)$ we have, $\nu-a.s.$,  $$b^\nu = \int_0^. E_{\nu}\left[\xi^V_t \middle| \mathcal{F}_t^\nu\right]  dt $$ where $(\xi^V_t)_{t\in[0,1)}$ is the stochastic process on $(W,\mathcal{B}(W)^\nu,\nu)$  defined by $$(\xi^V_t)_{t\in[0,1)} : (t,\omega)\in [0,1) \times W \to \xi^V_t(\omega) := \frac{\omega_1-\omega_t}{1-t} + \int_t^1 \frac{(1-s)}{(1-t)} \nabla V(\omega_s) ds \in \mathbb{R}^d$$  \end{proposition}\nproof By the stochastic Euler-Lagrange condition the process $(A_t)$ defined by $$A_t(\omega):= v_t^\nu + \int_0^t (\nabla V)(W_s) ds $$ is a $(\mathcal{F}_t^\nu)-$martingale so that we have $\nu\otimes \lambda-a.s.$\begin{eqnarray*}A_t(\omega) & = & \frac{1}{1-t}E_\nu\left[\int_t^1 A_\sigma d\sigma \middle|\mathcal{F}_t^\nu\right] \\ & = & E_\nu\left[\frac{W_1-W_t}{1-t}\middle|\mathcal{F}_t^\nu\right] + \frac{1}{1-t}E_\nu\left[\int_t^1 \left(\int_0^s (\nabla V)(W_\sigma) d\sigma \right) ds \middle|\mathcal{F}_t^\nu \right]\end{eqnarray*}where the last line is obtained by noticing that, from the  definition of $v^\nu$, we have $$ E_\nu\left[\frac{W_1-W_t}{1-t}\middle|\mathcal{F}_t^\nu\right] = E_\nu\left[\frac{\int_t^1 v_s^\nu ds}{1-t}\middle|\mathcal{F}_t^\nu\right]$$  Hence we obtain  $\nu\otimes \lambda-a.s.$ \begin{equation} \label{mpls} v_t^\nu =E_\nu\left[\frac{W_1-W_t}{1-t}\middle|\mathcal{F}_t^\nu\right] -\int_0^t (\nabla V)(W_s) ds + \frac{1}{1-t}E_\nu\left[\int_t^1 \left(\int_0^s (\nabla V)(W_\sigma) d\sigma \right) ds \middle|\mathcal{F}_t^\nu \right] \end{equation} and the result follows by integrating by parts. \nqed

\subsection{Critical processes of classical Lagrangians and  systems of stochastic differential equations}
\begin{theorem}Let $\widehat{\eta}\in \mathcal{P}_{\mathbb{R}^d}$ be a Borel probability on $\mathbb{R}^d$. Assume that $(X, Y)$ satisfies the system\begin{eqnarray}\label{1fbs} dX_t  & = & \sigma_t(X) dB_t + Y_t dt  \\
\label{2fbs}  dY_t & = & dZ_t - (\nabla V)(t,X_t) dt \\
\label{3fbs} Law(X_0) & = &  \widehat{\eta} 
\end{eqnarray}
on some complete stochastic basis $(\Omega,\mathcal{A}, (\mathcal{A}_t)_{t\in[0,1]}, \mathcal{P})$, where $(Z_t)$ is a \textit{c\`{a}dl\`{a}g} $(\mathcal{A}_t)_{t\in[0,1)}-$ martingale, $(B_t)$ an $(\mathcal{A}_t)-$Brownian motion, and $(X_t)$ is $(\mathcal{A}_t)-$adapted. Then $$X_\star\mathcal{P} \in G_V(\mathbb{R}^d)$$ and conversely if  $\nu \in G_V(\mathbb{R}^d)$  with ${W_0}_\star \nu = \widehat{\eta}$ then $(W_t, v^\nu_t)$ satisfy a system of this form with $(\sigma_t)$ such that $\alpha^\nu_. = \sigma_. \sigma^{\dagger}_.$ on the space $(W,\mathcal{B}(W)^\nu,\nu)$ with the filtration $(\mathcal{F}_t^\nu)$ or on one of its extensions (see \cite{I-W}). \end{theorem}
\nproof
We assume that $(X,Y, Z)$ is a solution to~(\ref{1fbs})~(\ref{2fbs})~(\ref{3fbs}) on a probability space $(\Omega,\mathcal{A},\mathcal{P})$ with a filtration $(\mathcal{A}_t)_{t\in[0,1]}$ for a $(\mathcal{A}_t)-$ Brownian motion $B$. We set $$\nu:=X_\star\mathcal{P}$$ By Proposition~\ref{lemmaforback} condition~(\ref{1fbs}) implies that $\nu \in \mathbb{S}$ with $$\int_0^. \alpha^\nu_t dt =  \int_0^. \sigma_t \sigma^{\dagger}_tdt$$ and  $\nu-a.s.$ $$b^\nu\circ X = \int_0^. E_{\mathcal{P}}\left[Y_t \middle|\mathcal{G}_t^X\right]dt $$ where $b^\nu$ denotes the finite variation part of $\nu$ (see Definition~\ref{defvelocity}), and where $\left(E_{\mathcal{P}}\left[Y_t \middle|\mathcal{G}_t^X\right]\right)$ denotes a \textit{c\`{a}dl\`{a}g} modification of the optional projection of $(Y_t)$ on the usual augmentation $(\mathcal{G}_t^X)$  of the filtration generated by $X$. We now take $(t,\omega)\in [0,1]\times W \to v_t^\nu\in \mathbb{R}^d$ to be \textit{c\`{a}dl\`{a}g}, $(\mathcal{F}_t^\nu)$-adapted and such that $\mathcal{P}-a.s.$ for all $t\in[0,1]$ $$v_t^\nu\circ X = E_{\mathcal{P}}\left[Y_t \middle|\mathcal{G}_t^X\right]$$ and we set \begin{equation} \label{ntnu} N_t^\nu := v_t^\nu + \int_0^t (\nabla V)(\sigma, W_\sigma) d
\sigma \end{equation} We want to prove that $(N_t^\nu)$ is a $(\mathcal{F}_t^\nu)-$martingale.  First note that by definition $(Z_t)$ is a $(\mathcal{A}_t)$ 
martingale and that for any $t\in[0,1]$ $\mathcal{G}_t^X\subset \mathcal{A}_t$. Therefore, for any $s\leq t$, $\mathcal{P}-a.s.$ \begin{equation}  \label{condmartfb}E_{\mathcal{P}}\left[ Z_t \middle| \mathcal{G}_s^X\right] = E_{\mathcal{P}}\left[ Z_s \middle| \mathcal{G}_s^X\right] \end{equation}  On the other hand by~(\ref{ntnu}) and~(\ref{veqlfb})  we have, for 
 $0\leq s \leq t <1$, $ \mathcal{P}-a.s.$ $$ E_\nu\left[N_t^\nu \middle| \mathcal{F}_s^\nu\right]\circ X   =   E_{\mathcal{P}}\left[N_t^\nu\circ X \middle| \mathcal{G}_s^X \right]  =    E_\mathcal{P}\left[ \left( E_\mathcal{P}[Y_t|\mathcal{G}_t^X]  + \int_0^t (\nabla V)(\sigma, X_\sigma) d\sigma\right)  \middle|\mathcal{G}_s^X\right]  =    E_\mathcal{P}\left[\widetilde{Z}_t  \middle|\mathcal{G}_s^X\right]  $$ where  $\widetilde{Z}_t:= Z_t -Z_0+Y_0$. 
Together with~(\ref{condmartfb}) and since $(N_t^\nu)$ is $(\mathcal{F}_t^\nu)$ adapted we derive, for $s\leq t$ $\mathcal{P}-a.s.$ $$E_\nu\left[N_t^\nu| \mathcal{F}_s^\nu\right]\circ X =  E_\mathcal{P}\left[\widetilde{Z}_t  \middle|\mathcal{G}_s^X\right] = E_\mathcal{P}\left[\widetilde{Z}_s \middle|\mathcal{G}_s^X\right]  = N^\nu_s\circ X$$so that $(N_t^\nu)$ is a $(\mathcal{F}_t^\nu)-$ martingale on $(W,\mathcal{B}(W)^\nu, \nu)$. The converse follows from the definition. Indeed  take $Z_t=N_t^\nu$, where $N_t^\nu$ is the martingale of the stochastic Euler-Lagrange condition. Then $$dv_t^\nu = dZ_t  +  (\nabla V)(t, W_t)dt$$ and on the other hand,  $\nu-a.s.$ $$dW_t = dM_t^\nu + v_t^\nu dt$$ so that the result follows by the martingale representation theorem (see \cite{I-W}) which can be applied on $(W, \mathcal{B}(W)^\nu, \nu)$ or on one of its extensions if $\sigma^\nu$ is degenerated.
\nqed
\begin{remarkk}
It is an interesting problem to determine conditions that assure the existence of solutions for a system  of type~(\ref{1fbs})-(\ref{3fbs}). Adding a final  
condition for the process $Y$ leads us to the study of forward-backward stochastic differential  systems (we refer, for  example to \cite{[n]}).
\end{remarkk}

The next result shows the existence of a probability   satisfying the Euler-Lagrange condition for a classical Lagrangian with force $V=0$.
\begin{proposition}
Let  $\gamma$ be a Borel probability on $\mathcal{P}_{\mathbb{R}^d\times \mathbb{R}^d}$ whose first marginal $\pi_\star \gamma$ is denoted by $\nu_0\in \mathcal{P}_{\mathbb{R}^d}$, and let $\mu_{\nu_0}:=\int_{\mathbb{R}^d} \nu_0(dx)\mu^{x}$ where $\mu^x$ is the Wiener measure starting from $x\in \mathbb{R}^d$. Assume also that \begin{equation} \label{leo} l(\gamma):= \inf( \mathcal{H}(\nu | \mu_{\nu_0}) : \nu \in \mathcal{P}_W,  ({W_0}\times W_1)_\star \nu = \gamma) <\infty  \end{equation} where $$\mathcal{H}(\nu | \mu_{\nu_0}):= E_\nu\left[\ln \frac{d\nu}{d\mu_{\nu_0}}\right]$$ denotes the relative entropy. Then there exists a unique probability $\nu^\star$ which attains the infimum of~(\ref{leo}) and $\nu^\star \in G_0(\mathbb{R}^d)$.
\end{proposition}
\nproof
By a classical application of the Dunford$-$Pettis theorem, the relative entropy with respect to $\nu_0$ has compact level sets (for the weak convergence in measure). Moreover it is strictly convex.  Since $\{\nu \in \mathcal{P}_W  : 
({W_0}\times W_1)_\star \nu \}$ is closed (for the weak convergence in measure) and convex, as soon as $l
(\gamma) < \infty$ the infimum is attained by a unique probability $\nu^\star$. If $\mathcal{H}(\nu | \mu_{\nu_0})<
\infty$ then in particular $\nu <<\mu_{\nu_0}$ and by the Gisanov theorem $\nu \in \mathbb{S}_B$ (see Corollary~
\ref{corollarycontrole}). On the other hand by the celebrated formula of \cite{F3} if $\mathcal{H}(\nu | \mu_{\nu_0})<\infty$ we have $$ 2\mathcal{H}(\nu | \mu_
{\nu_0}) = E_\nu\left[\int_0^1 |v_s^\nu|^2 ds\right]$$   Thus we obtain $$ l(\gamma)= \inf\left(\left\{E_\nu\left[\int_0^1 \frac{|v_s^\nu|^2}{2} ds\right] , \nu \in  \mathbb{S}_B \cap \mathbb{S}_{\gamma}\right\}\right)$$ and the result follows by Corollary~\ref{corollarycontrole}.  \nqed

\begin{examples}
\label{fin}
Let us mention the following examples and counterexamples which follow by simple calculus :
\begin{enumerate}[(i)]
\item For any $x,y\in \mathbb{R}^d$, $\mu_{x,y}\in G_0(\mathbb{R}^d)$, where $\mu_{x,y}$ denotes the law of the pinned Brownian motion such that $W_0=x$ and $W_1=y$.

\item For $d=1$, let $\alpha \in[0,1)$ and $\mu^0$ be the standard Wiener measure with $\mu^0-$a.s. $W_0=0$. Define $\nu$ to be the probability which is absolutely continuous with respect to $\mu^0$ with density  given by $$\frac{d\nu}{d\mu^0}:= \frac{(W_1 - W_\alpha)^2}{1-\alpha}$$ Then Clark$-$Ocone formula of Malliavin calculus shows that $\nu \in G_0(\mathbb{R}^d)$ and that it is the law of the non-Markovian stochastic differential equation $$dX_t = dB_t +  2  1_{[\alpha,1]}(s) \frac{X_s-X_\alpha}{1-s +(X_s-X_\alpha)^2} ; X_0=0$$ Now if we denote by $\nu^\star$ the probability defined by $$\frac{d\nu^\star}{d\mu^0}= \frac{d{W_1}_\star\nu }{d{W_1}_\star \mu^0}(W_1)$$ and since $\mu^0-a.s.$ $W_0=0$, by Jensen's inequality $\nu^\star$ solves $$\inf( \mathcal{H}(\eta |\mu_0): (W_0\times W_1)_\star \eta =(W_0\times W_1)_\star\nu )$$ By the strict convexity of the entropy we obtain $$ \mathcal{H}(\nu^\star |\mu_0)<  \mathcal{H}(\nu|\mu_0)  <\infty$$  This proves that, even for the Lagrangian $\mathcal{L}^0= \frac{|v|^2}{2}$, there may exist several elements of $G_0(\mathbb{R})$ of finite entropy with the same joint laws for $W_0$ and $W_1$.

\item Assume that $$u: (t,x) \in [0,1]\times \mathbb{R}^d \to \mathbb{R}^d$$ is a $C^{1,2}$ function, which is essentially bounded with an essentially bounded gradient, and that satisfies the incompressible Navier-Stokes equation
\begin{equation} \label{NSeq} \partial_t u + (u.\nabla)u = -\nabla p + \frac{\Delta u}{2} ,~~~~~~~ \hbox{div}~ u=0\end{equation} By setting $$\frac{d\nu}{d\mu^0} :=  \exp\left(- \int_0^1u(1-t,W_t) dW_t- \frac{1}{2} \int_0^1 |u(1-t, W_t)|_{\mathbb{R}^d}^2 dt \right)$$ we obtain a probability which is equivalent to the Wiener measure $\mu^0$. In particular, by the Girsanov theorem and Paul Levy's criterion, we have that $$M_t^\nu := W_t-W_0 + \int_0^t u(1-s, W_s) ds$$ is a $(\mathcal{F}_t^\nu)-$Brownian motion on $(W,\mathcal{B}(W)^\nu,\nu)$ so that $\nu\in \mathbb{S}$ and $\lambda\otimes \nu-a.s.$ $$v_t^\nu = -u(1-t, W_t)$$ (see Definition~\ref{defvelocity}). Since $u$ satisfies~(\ref{NSeq}), the It\^{o} formula directly implies that $\nu \in G_V(\mathbb{R}^d)$ for $$V= p(1-t,x)$$  i.e. $$t\to v_t^\nu +\int_0^t \nabla V(s,W_s) ds $$ is a $(\mathcal{F}_t^\nu)-$martingale.

\end{enumerate}
\end{examples}


\begin{thebibliography}{9}


\bibitem{AM} \textsc{Abraham, R.}, \textsc{Marsden, J.E.}, \textit{Foundations of Mechanics}, Addison-Wesley Publishing Company; 2nd edition (1980)
\bibitem{ARNOLD6} \textsc{Arnold,  V.I.}, \textit{Mathematical methods of classical mechanics},  Graduate Texts in Mathematics 60, Springer-Verlag; 2nd edition (1989)
\bibitem{DM2} \textsc{Dellacherie, C.}, \textsc{Meyer, P.A.}, \textit{Th\'{e}orie des martingales - Chapitres 5, 8},  Hermann (????)
\bibitem{EMERY} \textsc{\'{E}mery, M.}, \textit{En cherchant une caract\'erisation variationnelle des martingales}, S\'{e}minaire de probabilit\'{e}s de Strasbourg, 22 (1988), p. 147-154 


\bibitem{F3}
\textsc{F\"{o}llmer, H.}, \textit{Random fields and diffusion processes.},   Ecole d' \'{e}t\'{e} de Saint Flour XV--XVII  (1988)  


\bibitem{I-W} \textsc{Ikeda, N.}, and \textsc{Watanabe, S.}, 
\textit{Stochastic Differential Equations and Diffusion Processes}, North
Holland, Amsterdam (Kodansha Ltd., Tokyo) (1981)

\bibitem{JACOD}
\textsc{Jacod, J.}, \textsc{Shiryaev, A.N.}, \textit{Limit theorems for stochastic processes}, Springer Verlag (1987)



\bibitem{LANDAU}

\textsc{Landau, L.},\textsc{Lifchitz}, \textit{Physique Th\'{e}orique 1 M\'{e}canique}, Editions Mir Moscou U.R.S.S.; 4th edition (1988)


\bibitem{ABSAPI}
\textsc{Lassalle, R.}, \textit{Invertibility of adapted perturbations of the identity on abstract Wiener space},  J. Func. Anal. 262(6):2734-2776 (2012)


\bibitem{RL-locinv}
\textsc{Lassalle, R.},\textit{Local invertibility of adapted shifts on Wiener space, under finite energy condition}, Stochastics and Stochastics Reports (2012) 

\bibitem{RLASU} \textsc{Lassalle, R.}, \textsc{\"{U}st\"{u}nel, A.S. } \textit{Local invertibility of adapted shifts on Wiener space, and related topics},  Springer Proceedings in Mathematics \& Statistics vol. 34 (2013)

\bibitem{RLZ} \textsc{Lassalle, R.}, \textsc{Zambrini, J.C.}, \textit{A weak approach to the stochastic deformation of classical mechanics}, submitted paper (2014)



\bibitem{LEO} \textsc{Leonard, C.}, \textit{A survey of the Schr\"{o}dinger problem and some of its connections with optimal transport}, Discrete and Cont. Dyn. systems A, 34(4), 1533 (2014)
\bibitem{leo} \textsc{Leonard, C. },  \textsc{Roelly, S.},  \textsc{Zambrini, J.C.},  \textit{Reciprocal processes. A measure-theoretical point of view}, to be published in Probability Surveys.
\bibitem{[n]}
\textsc{Ma, J.}, \textsc{Yong. J.}, \textit{Forward-Backward Stochastic Differential  Equations and their Applications}, Lecture Notes in Mathematics 1702, Springer (2007)


\bibitem{MAL2}
\textsc{Malliavin, P.}, \textit{Int\'{e}gration et probabilit\'{e}s}, MASSON (1982) 

\bibitem{Meyer-Zheng}
\textsc{Meyer, P.A.}, \textsc{Zheng, W.A.},  \textit{Tightness criteria for laws of semi-martingales}, Ann. Inst. Poincar\'{e}, (1984)



\bibitem{Mika1}
\textsc{Mikami ,T.},\textsc{Thieullen,M.}, \textit{Duality Theorem for Stochastic Optimal Control Problem}, Stoch. Proc. Appl. 116,1815-1835, MR230760  (2006) 

\bibitem{Mika3}
\textsc{Mikami, T.}, \textsc{Thieullen, M.}, \textit{Optimal transportation problem by stochastic optimal control} (2005) 

\bibitem{Mika4} \textsc{Mikami, T.},    \textit{Optimal Transportation Problem as Stochastic Mechanics}  (2008) 


\bibitem{RR}
\textsc{Robertson, A. P.},  \textsc{Robertson, W.}, \textit{Topological Vector Spaces}, Cambridge University Press, (1964)

\bibitem{TOUZI}
\textsc{Tan, X.}, \textsc{Touzi, N.}, \textit{Optimal transportation under controlled stochastic dynamics}, The Annals of Probability (2013)
 \bibitem{ZT}
 \textsc{Thieullen, M.} \textsc{Zambrini, J.C.}  \textit{Probability and quantum symmetries I, the theorem of Noether in Schr\"{o}dinger's euclidean quantum mechanics}, Ann. Inst. H.Poincar\'e, Phys. theo. Vol 67, N 3, 1997, p.297

\bibitem{Tsi}
\textsc{Tsirelson , B.S.}, \textit{An example of stochastic differential equation having no strong solution}. Theor. Prob. Appl. {\bf 20}, (1975) , p. 416--418


 \bibitem{ASU-3}
\textsc{ \"Ust\"unel, A. S.}, 
\textit{  Entropy, invertibility and variational calculus
 of adapted shifts on Wiener space}, J. Funct. Anal. 257,  no. 11, 3655--3689. (2009)





\bibitem{Zamb1}
\textsc{Zambrini, J.C.}, \textit{Stochastic mechanics according to E. Schr\"{o}dinger}, Phys. Rev. A, volume 33, number 3.(1986) 

\bibitem{Zamb3}
\textsc{Zambrini, J.C.},  \textit{Variational processes and stochastic versions of mechanics}, J. Math. Phys. 27 (9). (1986)

\bibitem{ZHENG}
\textsc{Zheng, W.A.}, \textit{Tightness results for laws of diffusion processes application to stochastic mechanics}, Annales de l'I.H.P., section B, tome 21, (1985)

\end{thebibliography}
\end{document}